\documentclass[12pt]{article}
\usepackage{amssymb,amsmath}
\usepackage[latin1]{inputenc}
\numberwithin{equation}{section}
\numberwithin{figure}{section}

\def\ba{\begin{array}}
\def\ea{\end{array}}

\def\Re{\mathop{\mathrm{Re}}}
\def\Im{\mathop{\mathrm{Im}}}
\def\ts{\textstyle}

\def\pl{\partial}

\def\set#1#2{ \{ #1 \:|\: #2 \} }
\def\bigset#1#2{ \bigl\{ #1 \:\big|\: #2 \bigr\} }
\def\Bigset#1#2{ \Bigl\{ #1 \:\Big|\: #2 \Bigr\} }

\def\i{\mathrm i}
\def\e{\mathrm e}
\def\d{\mathrm d}
\def\dd{\;\! \d}
\def\eps{\varepsilon}
\def\loc{\mathrm{loc}}
\def\supp{\mathrm{support}}
\def\ae{\text{a.a.\ }}

\def\R{\mathbb R} \def\C{\mathbb C} \def\Z{\mathbb Z}
\def\N{\mathbb N}

 \def\cD{\mathcal D}
 \def\cF{\mathcal F}

 \def\cL{\mathcal L}
 \def\cN{\mathcal N}
\def\cO{\mathcal O} 
 
\def\cS{\mathcal S}

\def\bbC{\mathbb C} 
 
 \def\bbH{\mathbb H}
 
 \def\bbL{\mathbb L}
 \def\bbN{\mathbb N}
 \def\bbP{\mathbb P}
\def\bbQ{\mathbb Q} \def\bbR{\mathbb R}

\def\wh{\widehat}
\def\wt{\widetilde}
\def\euler{\gamma_{\text{\sc E}}}

\def\build#1_#2^#3{\mathrel{
  \mathop{\kern 0pt#1}\limits_{#2}^{#3}}}

\newcommand{\QED}{\mbox{}\hfill 
\raisebox{-0.2pt}{\rule{5.6pt}{6pt}\rule{0pt}{0pt}}\medskip\par} 
\newtheorem{lemma}{Lemma}[section]
\newtheorem{proposition}[lemma]{Proposition}
\newtheorem{theorem}[lemma]{Theorem}
\newenvironment{example}{{\noindent \bf Example.}}{\smallskip}
\newenvironment{remark}{{\noindent \bf Remark.}}{\smallskip}
\newenvironment{remarks}{{\noindent \bf Remarks.}}{\smallskip}
\newenvironment{definition}{{\noindent \bf Definition.}}{\smallskip}
\newenvironment{proof}{{\it \noindent Proof. }}{\QED}
\newcommand{\Ref}[1]{\mbox{\rm (\ref{#1})}}            

\newdimen\texpscorrection
\newdimen\figcenter
\def\figurewithtex#1#2#3#4#5#6{\null
  {\goodbreak\figcenter=\hsize\relax
  \advance\figcenter by -#4truecm
  \divide\figcenter by 2
  \begin{figure}[hbt]
  \vskip #3truecm\noindent\hskip\figcenter
  \includegraphics{#1}{\hskip\texpscorrection\input #2 }
  \vskip 0.0truecm\noindent \caption{{\footnotesize #6}}\label{#5}
  \end{figure}}}

\def\point#1 #2 #3 {\rlap{\kern #1 truecm \raise #2 truecm \hbox{#3}}}

\def\epsilon{\eps}
\def\FOR{\quad\text{for }}
\def\AND{\text{ and }}
\def\Eone{{\mathrm E}_1}

\begin{document}
\title{Convergence results for a coarsening model\\
using global linearization}
\author{Thierry Gallay\\
\normalsize Institut Fourier\\
\normalsize Universit\'e de Grenoble I\\
\normalsize BP 74\\
\normalsize F-38402 Saint-Martin d'H\`eres\\
\and
Alexander Mielke\\
\normalsize Mathematisches Institut A\\
\normalsize Universit\"at Stuttgart \\
\normalsize Pfaffenwaldring 57\\
\normalsize D-70569 Stuttgart}
\date{September 13, 2002}
\maketitle

\begin{abstract}
We study a coarsening model describing the dynamics of interfaces
in the one-dimensional Allen-Cahn equation. Given a partition of
the real line into intervals of length greater than one, the model
consists in constantly eliminating the shortest interval of
the partition by merging it with its two neighbors. We show that
the mean-field equation for the time-dependent distribution of
interval lengths can be explicitly solved using a global linearization
transformation. This allows us to derive rigorous results on the
long-time asymptotics of the solutions. If the average length of
the intervals is finite, we prove that all distributions approach 
a uniquely determined self-similar solution. We also obtain global 
stability results for the family of self-similar profiles which
correspond to distributions with infinite expectation. 
\end{abstract}


\section{Introduction}

Consider a domain $D \subset \bbR^n$ which is divided into a large
number of subdomains (or cells) of different sizes, separated by
domain walls, and assume that the system evolves in such a way that
the larger subdomains grow with time while the smaller ones shrink
and eventually disappear. In particular, the average size of the
cells increases, so that the subdivision of $D$ becomes rougher
and rougher. Such a {\it coarsening} dynamics is observed in many 
physical situations, especially near a phase transition when a system 
is quenched from a homogeneous state into a state of coexisting phases.
Typical examples are the formation of microstructure in alloy 
solidification \cite{LiS61} and the phase separation in lattice spin 
systems \cite{De97,KBN97}. Closely related to coarsening is the 
coagulation (or aggregation) process which describes the
dynamics of growing and coalescing droplets \cite{DGY91,PeR92,Vo85}.
In this case, the system consists of a large number of particles 
of different masses which interact by forming clusters. Again, 
the total mass is preserved, so that the average mass per cluster 
increases with time. 

Given a coarsening or a coagulation model, the main task is to 
predict the long-time evolution of the size distribution of the
cells, or the mass distribution of the clusters. In many cases, 
experiments and numerical calculations show that this behavior 
is asymptotically self-similar: the system can be described by 
a single length scale $\cL(t)$, and the distribution approaches 
the scaling form $\cL(t)^{-1} \Phi(x/\cL(t))$ as $t \to \infty$. 
The profile $\Phi$ and the asymptotics of $\cL(t)$ can sometimes 
be determined exactly \cite{NaK86,BDG94}. However, even in simple 
situations, it is very difficult to prove that the size 
distribution actually converges to a self-similar profile. 

In this work, we consider a simple coarsening model related to the 
one-dimensional Allen-Cahn equation $\pl_t u = \pl_x^2 u + 
\frac12(u -u^3)$, where $x \in \R$. The equilibria of this system 
are the homogeneous steady states $u = \pm 1$, together with the
kinks $u(x) = \pm\tanh(x/2)$ which represent domain walls separating 
regions of different ``phases''. If $u$ is any bounded solution of 
this equation, then for $t > 0$ sufficiently large the graph of
$u(t,\cdot)$ will typically look like a (countable) family of kinks 
separated by large intervals on which $u \approx \pm 1$. If we 
denote by $x_j(t)$ the position of the $j^{\rm th}$ kink and
if we assume that $x_{j+1}(t) - x_j(t) \gg 1$ for all $j \in \Z$, 
a rigorous asymptotic analysis shows that $\dot x_j \approx
F(x_{j+1}-x_j)-F(x_j-x_{j-1})$, where $F(y) = 24 \e^{-y}$ 
\cite{CaP89}. In other words, the positions of the domain walls
behave like a system of point particles with short range attractive 
pair interactions. Thus, on an appropriate time scale, only the
closest pairs of kinks will really move; in such pairs, kinks will 
attract each other until they eventually annihilate.

This kink dynamics suggests the following coarsening model
\cite{NaK86,DGY91,CaP92,BDG94,RuB94,BrD95,CaP00}. Consider a 
partition of the real line $\R$ into a countable union of 
disjoint intervals $I_j$, with $\ell(I_j) \ge 1$ for all $j \in \Z$.
In the previous picture, the intervals $I_j$ correspond to regions 
where $u$ is close to $\pm 1$. A dynamics on this configuration 
space is defined by iterating the following coarsening step: 
choose the ``smallest'' interval in the partition, and merge it
with its two nearest neighbors. This model clearly mimics the 
dynamics of the domain walls in the one-dimensional Allen-Cahn 
equation. However, proving that the formal procedure described
above actually defines a well-posed evolution (e.g. for almost 
all initial configurations) and investigating its statistical
properties after many coarsening iterations is a non-trivial task, 
which has not been accomplished so far. Instead, the coarsening
model has been studied in the {\it mean field} approximation, 
which consists in merging the minimal interval not with its true 
neighbors, but with two intervals chosen at random in the
configuration $\{I_j\}_{j\in \Z}$. This approximation is valid
provided the lengths of consecutive intervals stay uncorrelated 
during the coarsening process, an assumption that seems reasonable
\cite{BDG94}. 

Under this hypothesis, it is possible to write a closed evolution 
equation for the distribution $f(t,x)$ (per unit length) 
of intervals of length $x \ge 1$
at time $t$ \cite{CaP92}. Denoting by $N(t) = \int_0^\infty 
f(t,x)\dd x$ the total number of intervals per unit length, 
and by $\cL(t)$ the 
length of the smallest interval, the equation reads
\begin{equation}\label{eq:feq}
  \pl_t f(t,x) = \frac{\dot\cL(t)f(t,\cL(t))}{N(t)^2}
\left(\int_0^{x-\cL(t)} f(t,y)
  f(t,x{-}y{-}\cL(t))\dd y - 2f(t,x)N(t)\right),
\end{equation}
for $x \ge \cL(t)$, whereas $f(t,x) = 0$ for $x < \cL(t)$ by 
the definition of $\cL(t)$. By construction, $N(t)$ decreases with 
time, while the total length of the intervals $\int_0^\infty 
xf(t,x)\dd x$ is conserved. 
\medskip

We prefer to work with the distribution density $\rho(t,x) =
f(t,x)/N(t)$, which satisfies $\rho(t,x) = 0$ for $x < \cL(t)$ and the 
normalization $\int_0^\infty \rho(t,x)\dd x = 1$ for all $t$. The 
evolution equation for $\rho$ reads
\begin{equation}\label{eq:roeq_new}
  \pl_t \rho(t,x) = \dot{\cL}(t)\rho(t,\cL(t)) \int_0^{x-\cL(t)} \rho(t,y)
  \rho(t,x{-}y{-}\cL(t))\dd y  \FOR x \ge \cL(t).
\end{equation}
Of course, systems \Ref{eq:feq} and \Ref{eq:roeq_new} are equivalent. 
In particular, once the density $\rho(t,x)$ is known, the total number 
$N(t)$ can be recovered by solving the ordinary differential equation 
$\dot N(t) = -2\dot{\cL}(t)\rho(t,\cL(t))N(t)$, and the distribution 
$f(t,x)$ is then given by $N(t)\rho(t,x)$. 

It is important to note that equations~\Ref{eq:feq}, \Ref{eq:roeq_new}
are invariant under reparametrizations of time. As a consequence, 
the minimal length $\cL(t)$ is not determined by the initial 
data, but can be prescribed to be an arbitrary (increasing) function 
of time. In \cite{CaP92}, the authors define an ``intrinsic time'' by 
imposing the relation $f(t,\cL(t))\dot\cL(t) = 1$, which means that 
the number of merging events per unit time is constant. We find it 
more convenient to use the ``coarsening time'' defined by the 
simple relation $\cL(t) = t$. In other words, we choose to 
parameterize the coarsening process by the length of the smallest 
remaining interval, forgetting about how much physical time elapses 
between or during the merging events. With our choice,  
equation~\Ref{eq:roeq_new} becomes
\begin{equation}\label{eq:roeq}
  \pl_t \rho(t,x) = \rho(t,t) \int_0^{x-t} \rho(t,y)
  \rho(t,x{-}y{-}t)\dd y  \FOR x \ge t.
\end{equation}
Since we do not allow for intervals of length smaller than $1$, 
we impose our initial condition at time $t=1$: $\rho(1,x) =
\rho_1(x)$. 
 
The aim of this paper is to show that the dynamics of \Ref{eq:roeq}
can be completely understood using a global linearization 
transformation. As a consequence, we are able to prove that 
solutions of \Ref{eq:roeq} satisfying $\int_0^\infty x\rho(t,x)
\dd x < \infty$ approach a non-trivial self-similar profile as 
$t \to \infty$. To achieve this goal, we first rewrite \Ref{eq:roeq} 
in similarity coordinates by setting
\[
  \rho(t,x) = \frac{1}{t} \eta(\log t, x/t), \quad\text{or} \quad 
  \eta(\tau,y) = \e^\tau\rho(\e^\tau,\e^\tau y),
\]
where $\tau = \log t \ge 0$ and $y = x/t \in [1,\infty)$. 
Then the rescaled density $\eta(\tau,\cdot)$ lies in the 
time-independent space
\begin{equation}\label{eq:Pdef}
  \bbP = \Bigset{\eta \in L^1((1,\infty),\R_+)}{
  \int_1^\infty \eta(y)\dd y = 1}, 
\end{equation}
which is a closed convex subset of $L^1((1,\infty))$. Moreover, 
\Ref{eq:roeq} is transformed into the autonomous evolution 
equation
\begin{equation}\label{eq:etaeq}
  \pl_\tau \eta(\tau,y) = \pl_y\big(y \, \eta(\tau,y)\big) + 
  \eta(\tau,1) \int_1^{y-2} \eta(\tau,z)\eta(\tau,y{-}z{-}1)\dd z 
  \FOR y \ge 1. 
\end{equation}
In Section \ref{s:S} we show that, for all initial data 
$\eta_0 \in \bbP$, \Ref{eq:etaeq} has a unique global solution 
$\eta \in C^0([0,\infty),\bbP)$ with $\eta(0) = \eta_0$.

We now define a nonlinear map $\cN : \bbP \to L^1_\loc([1,\infty),\R_+)$ 
by
\[
  \cN = \cF^{-1} \circ \phi \circ \cF,
\]
where $\cF$ is the Fourier transform and $\phi(z) = \frac12 \log
\frac{1+z}{1-z}$. If $\eta(\tau,\cdot)$ is a solution of \Ref{eq:etaeq} 
in $\bbP$, a direct calculation reveals that $w(\tau,\cdot) = 
\cN(\eta(\tau,\cdot))$ satisfies the linear equation 
$\pl_\tau w(\tau,y) = \pl_y (y w(\tau,y))$. As a consequence, 
\[
  w(\tau,y) = (S_\tau w_0)(y) = \left\{\ba{cc} \e^\tau w_0(\e^\tau y)
  &\text{if }y \geq 1,\\ 0&\text{if }y<1,\ea\right.
\]
where $w_0 = \cN(\eta_0)$. It follows that any solution 
$\eta \in C^0([0,\infty),\bbP)$ of \Ref{eq:etaeq} satisfies 
$\cN(\eta(\tau)) = S_\tau \cN(\eta_0)$ for all $\tau \ge 0$. 
In other words, the nonlinear evolution defined by \Ref{eq:etaeq}
is {\it conjugated} (via the map $\cN$) to the linear semigroup 
$(S_\tau)$. Thus, the difficulty of solving \Ref{eq:etaeq} is 
carried over to the study of the mapping $\cN$ and of its inverse 
$\cN^{-1} = \cF^{-1} \circ \phi^{-1} \circ \cF$. Although the 
properties of these maps are not fully understood, it possible to 
obtain some information on them using the analyticity properties of 
the Fourier-Laplace transform. 

In Section \ref{SS} we investigate the steady states of \Ref{eq:etaeq}, 
which form a one-parameter family $\{\eta^*_\theta\}_{\theta \in \R}$. 
Here $\eta^*_\theta = \cN^{-1} (\frac{\theta}{2}w^*)$, where 
$w^*(y) = y^{-1}\mathbf 1_{\{y\ge 1\}}$. More explicitly, we have
\begin{equation}\label{eq:sssol}
 \wh{\eta^*_\theta}(\xi) = (\cF \eta^*_\theta)(\xi) = 
 \tanh\Bigl(\frac{\theta}{2}\, \Eone(\i \xi)\Bigr) \FOR \xi \in \R,
\end{equation}
where $\Eone$ is the exponential integral \cite{AS72}. 
We prove that $\eta^*_\theta \in \bbP$ if and only if $\theta\in
(0,1]$. Moreover, $\eta^*_1(y)$ decays exponentially as $y \to
\infty$, while $\eta^*_\theta(y) \sim y^{-(1+\theta)}$ if 
$0 < \theta < 1$. In particular, $\eta^*_1$ is the only steady state 
for which the average length $\int_1^\infty y\eta^*_1(y)\dd y$ is
finite. 

Finally, Section \ref{GC} is devoted to the convergence results. 
If the initial data $\eta_0 \in \bbP$ satisfy $y^\gamma \eta_0 \in 
L^2((1,\infty))$ for some $\gamma > 3/2$ (so that $\int_1^\infty 
y\eta_0(y)\dd y < \infty$), we prove that the corresponding solution 
of \Ref{eq:etaeq} converges exponentially to the steady state 
$\eta^*_1$:
\[
  \|y^{\gamma-1} (\eta(\tau)-\eta^*_1)\|_{L^2((1,\infty))} = 
  \cO(\e^{-(\gamma-3/2)\tau}) \FOR \tau\to \infty.
\]
In terms of the original variables, this shows that the density
$\rho(t,x)$ asymptotically approaches the self-similar solution 
$t^{-1}\eta^*_1(x/t)$ of \Ref{eq:roeq}. Moreover, the remainder is 
$\cO(t^{-(\gamma-3/2)})$, so that the convergence is very fast if 
$\gamma$ is large, i.e., the initial data decay rapidly at infinity. 
Similarly, if $0 < \theta < 1$ and if $\eta_0 \in \bbP$ satisfies 
$y^\gamma(\eta_0 {-} \nu\eta^*_\theta) \in L^2((1,\infty))$ for some 
$\gamma > \theta+ 1/2$ and some $\nu > 0$, we prove that the solution 
of \Ref{eq:etaeq} with initial data $\eta_0$ converges to the steady 
state $\eta^*_\theta$. 

To conclude this section, we briefly comment on previous results
and possible generalizations. The mean field equations \Ref{eq:feq}
and especially the self-similar solutions \Ref{eq:sssol} can 
be found in many physics papers \cite{NaK86,DGY91,BDG94,RuB94,BrD95}.
The first mathematical work is \cite{CaP92}, where the authors 
prove the existence of global solutions to \Ref{eq:feq}. They 
also show that the profile $\eta^*_1$ is a positive function 
(a crucial property that is tacitly assumed in the physics 
literature!) and study its asymptotic behavior as $y \to \infty$.
Our main contribution is the introduction of the linearization
transformation $\cN$ which allows to prove the convergence results. 
We also extend the analysis of \cite{CaP92} to the equilibria
$\eta^*_\theta$ with $0 < \theta < 1$.

The ``two-sided'' coarsening model discussed in this introduction
is clearly not the most general system to which our analysis applies.
For instance, we can consider the ``one-sided'' variant in which 
the minimal interval is only merged with {\it one} of its neighbors
\cite{CaP00}. More generally, we can assume that, for $j = 1,\dots,N$, 
the minimal interval has a probability $p_j$ of being merged with 
$j$ of its neighbors, where $p_1 + \dots + p_N = 1$. In the mean 
field approximation, this leads to an evolution equation similar
to \Ref{eq:roeq}, where the quadratic convolution in the right-hand side 
is replaced by a more general convolution polynomial. Except for a 
modified definition of the mapping $\cN$, this extension does not 
affect our analysis in any essential way. Therefore, in the rest of 
this paper, all results will be stated and proved in this
general situation. 

\paragraph*{Acknowledgments.} The authors are grateful for 
financial support through the French-German grant: PROCOPE 00307$\,$TK, {\it
  Attractors for Extended Systems}. 

\section{The coarsening equation and its solution}
\label{unscaled}

As is explained in the introduction, we shall study a general 
coarsening model for which the number of intervals involved in 
each merging event is not necessarily fixed. Instead, we allow 
for some randomness by choosing nonnegative real numbers 
$p_1, \dots, p_N$ satisfying $p_1 + \dots + p_N = 1$, where $p_j$
is interpreted as the probability for an interval of minimal length
to merge with $j$ other intervals. We define the polynomial
\[
   Q(z) = \sum_{j=1}^N p_j z^j, 
\]
which satisfies $Q(1) = 1$. The original coarsening model related 
to the Allen-Cahn equation corresponds to the particular case 
where $Q(z) = z^2$. 

If $\rho \in L^1(\R)$, we set
\begin{equation}\label{eq:nonlin}
  \bbQ[\rho] = \sum_{j=1}^N p_j \rho^{*j},
\end{equation}
where $\rho^{*j} = \rho * \rho * \dots * \rho$ ($j$ factors) and $*$ 
denotes the convolution product in $L^1(\R)$. In particular, we have
$\int_0^\infty \bbQ[\rho](x)\dd x = Q(\int_0^\infty \rho(x)\dd x)$. 
In what follows, we shall mainly use the space $\bbP$ of probability 
densities defined by \Ref{eq:Pdef}. Any $\rho \in \bbP$ can be 
extended to the whole real line by setting $\rho(x) = 0$ 
for $x < 1$. This natural extension, still denoted by $\rho$, will
be used in the sequel without further mention. As an example
of this abuse of notation, if $\rho \in \bbP$, we have $\bbQ[\rho] 
\in \bbP$ and $\supp(\bbQ[\rho]) \subset [2,\infty)$. 

The problem we are interested in can now be stated as follows. Given 
$\rho_1 \in \bbP$, find a density $\rho : [1,\infty)^2 \to \R_+$ 
satisfying $\rho(1,x) = \rho_1(x)$ for $x \ge 1$, $\rho(t,x) = 0$ for 
$1 \le x < t$, and
\begin{equation}\label{eq:U.1}
  \pl_t \rho(t,x) = \rho(t,t)\bbQ[\rho(t,\cdot)](x{-}t) \FOR x \ge t \ge 1. 
\end{equation}
If $Q(z) = z^2$, the evolution equation \Ref{eq:U.1} reduces to 
\Ref{eq:roeq}. 

By assumption, the density $\rho(t,x)$ is nonzero 
only in the sector $\set{(t,x) \in \R^2}{1 \le t \le x }$, where it 
satisfies \Ref{eq:U.1}. An important role will be played by the values 
of $\rho$ on the boundaries of this domain, namely the initial 
density $\rho_1$ and the trace of $\rho$ on the diagonal $x = t$, 
which we denote by $\alpha$:
\[
   \alpha(t) = \rho(t,t) \FOR  t \ge 1 .
\]
Any sufficiently smooth solution of \Ref{eq:U.1} satisfies 
$\rho(t,\cdot) \in \bbP$ for all $t \ge 1$ provided $\rho_1 \in 
\bbP$. Indeed, it is obvious from \Ref{eq:U.1} that $\rho$ stays
nonnegative. Moreover, if $m(t) = \int_t^\infty \rho(t,x)\dd x$, 
a direct calculation shows that
\begin{equation}\label{eq:momeq}
  \frac{\d}{\d t} m(t) = \alpha(t)\bigl(Q(m(t)) - 1\bigr) \FOR t \ge 1. 
\end{equation}
Therefore, if $m(1) = 1$, then $m(t) = 1$ for all $t \ge 1$. 

A very remarkable property of equation \Ref{eq:U.1} is that it can be 
explicitly solved using Fourier (or Laplace) transform. If $\rho 
\in \bbP$, we define
\[
  \wh\rho(\xi) = (\cF\rho)(\xi) = \int_1^\infty \e^{-\i \xi x} \rho(x)\dd x
  \FOR \xi \in \R .
\]
Then $\wh\rho \in C^0(\R,\C)$ satisfies $\wh\rho(0) = 1$, 
$|\wh\rho(\xi)| < 1$ for all $\xi \neq 0$, and $\wh\rho(\xi) \to 0$ 
as $\xi \to \pm\infty$. Moreover, $\wh\rho$ is a positive definite function 
(in the sense of Bochner). Since $\supp(\rho) \subset [1,\infty)$, 
the Fourier transform $\wh\rho$ can be continuously extended to 
the lower complex half plane
\[
  \bbL^- = \set{\xi\in \C }{ \Im\xi \leq 0 }. 
\] 
This extension (still denoted by $\wh\rho$) is analytic in the 
interior of $\bbL^-$ and satisfies the bound $|\wh\rho(\xi)| \le 
\e^{\Im\xi}$ for all $\xi \in \bbL^-$. 

\begin{remark}
The closely related Laplace transform is defined by
\[
  \wt\rho(p) = (\cL\rho)(p) = \int_1^\infty \e^{-p x} \rho(x)\dd x
  \FOR  \Re p \ge 0,  
\]
so that $\wt\rho(p)=\wh\rho(-\i p)$. In the sequel, we prefer 
using Fourier transform instead of Laplace because the inversion 
formula is more natural.  
\end{remark}

Applying Fourier transform to \Ref{eq:U.1} and using the fact 
that convolutions are turned into multiplications, we find the equation  
\begin{equation}\label{eq:U.2}
  \pl_t \wh\rho(t,\xi) = \alpha(t)\,\e^{-\i \xi t}\Big(Q(\wh\rho(t,\xi))
  -1\Big) \FOR  t \ge 1,
\end{equation}
where $\alpha(t)=\rho(t,t)$. To solve \Ref{eq:U.2}, we introduce the 
nonlinear complex transformation $\phi$ defined by 
\begin{equation}\label{eq:Uphi}
 \phi'(z) = \frac{1}{1-Q(z)}\:, \quad  \phi(0)=0.  
\end{equation}
Remark that $\phi'(z)=\sum_{k=0}^\infty [Q(z)]^k$, so that $\phi$ has 
a power series expansion with nonnegative coefficients whose radius of
convergence is equal to $1$. In particular, the map $\phi : [0,1) \to 
[0,\infty)$ is one-to-one and onto. Let $\psi = \phi^{-1}$ be the 
inverse map, which satisfies  
\begin{equation}\label{eq:Upsi}
  \psi'(w) = 1 - Q(\psi(w)), \quad  \psi(0)=0. 
\end{equation}
By construction, $\psi$ is analytic in a neighborhood of
the real positive axis. In the particular case where $Q(z) = z^2$, 
one finds 
\[
  \phi(z) = \frac12\,\log\frac{1+z}{1-z}\quad  \AND\quad   \psi(w) = \tanh(w). 
\]

Applying the nonlinear transformation $\phi$ simplifies equation
\Ref{eq:U.2} a lot. The function $\wh w(t,\xi) =\phi(\wh\rho(t,\xi))$, 
which is defined at least for $\Im \xi < 0$, satisfies the differential 
equation
\[
  \pl_t\wh w(t,\xi) = - \alpha(t)\e^{-\i \xi t} \FOR  t \ge 1, 
\]
which has the explicit solution
\begin{equation}\label{eq:Wexp}
  \wh w(t,\xi) = \wh w(1,\xi)-\int_{1}^t \alpha(s)\e^{-\i \xi s}\dd s \FOR 
  t \ge 1 \AND \Im \xi < 0. 
\end{equation}
Remark that $|\wh\rho(t,\xi)| \le \e^{t\Im \xi}$ for all $t \ge 1$ and
all $\xi \in \bbL^-$, because $\rho(t,\cdot) \in \bbP$ and $\supp
(\rho(t,\cdot)) \subset [t,\infty)$. Since $\phi(z) = z+\cO(|z|^2)$ as 
$z \to 0$, it follows that $|\wh w(t,\xi)| = |\phi(\wh\rho(t,\xi))|
\to 0$ as $t \to \infty$ if $\Im \xi < 0$. Thus, taking the limit 
$t\to \infty$ in \Ref{eq:Wexp}, 
we find $\wh w(1,\xi) = \int_1^\infty \alpha(t)\,\e^{-\i \xi t} \dd t$,
which in turn implies
\begin{equation}\label{eq:Wsol}
  \wh w(t,\xi) = \int_{t}^\infty \alpha(s)\e^{-\i \xi s}\dd s \FOR
   t \ge 1 \AND  \Im \xi < 0. 
\end{equation} 

This formula has a very nice interpretation. Let $\cN$ be the
nonlinear transformation defined (at least formally) by
\begin{equation}\label{eq:Ndef}
   \cN =  \cF^{-1} \circ \phi \circ \cF 
   \quad \text{or} \quad 
   \cN^{-1} = \cF^{-1} \circ \psi \circ \cF. 
\end{equation}
Setting $t = 1$ in \Ref{eq:Wsol}, we obtain $\phi(\wh\rho_1) = \wh\alpha$, 
that is $\alpha = \cN(\rho_1)$. In other words, the trace $\alpha(t) = 
\rho(t,t)$ is obtained from the initial density $\rho_1(x) = 
\rho(1,x)$ by applying the nonlinear map $\cN$. Moreover, if $U(t)$ 
is the linear operator defined for $t \ge 1$ by
\begin{equation}\label{eq:Udef}
  (U(t)w)(s) = \mathbf 1_{\{s\ge t\}} w(s) = \left\{\ba{cl}0&\text{if }s < t,\\
  w(s)&\text{if } s \ge t, \ea\right.
\end{equation}
then \Ref{eq:Wsol} reads $\wh w(t,\cdot)=\phi(\wh\rho(t,\cdot)) =
\cF(U(t)\alpha)$, which means $\cN (\rho(t,\cdot)) = U(t)\alpha$. Therefore,
the solution of \Ref{eq:U.1} satisfies
\begin{equation}\label{eq:U.3}
  \cN(\rho(t,\cdot)) = U(t)\cN(\rho_1) \FOR  t \ge 1.
\end{equation}
This shows that the dynamics of the nonlinear system \Ref{eq:U.1} is
conjugated via the nonlinear mapping $\cN$ to the linear evolution $U$.
Since $\cN(\rho_1)$ is the trace function defined by $\alpha(t) = 
\rho(t,t)$, it is very natural that the evolution of $\alpha$ is
obtained just by cutting off the history in $[1,t)$. 

It is not difficult to show that the map $\cN$ is well-defined 
on the space $\bbP$, cf.~\Ref{eq:Pdef}:

\begin{proposition}\label{S:Ndef}
If $\rho \in \bbP$, then $\cN(\rho) \in L^1_\loc([1,\infty),\R_+)$, 
and the mapping $\rho \mapsto \cN(\rho)$ is one-to-one. 
\end{proposition}

\begin{proof}
  For $\rho \in \bbP$ we construct $w=\cN(\rho)$ as follows.  Define $\wh w :
  \bbL^-_* \to \C$ by $ \wh w(\xi) = \phi(\wh\rho(\xi))$, where $\bbL^-_*=
  \bbL^- \setminus \{0\}$ .  We recall that $\wh\rho$ is continuous on
  $\bbL^-$, analytic in the interior of $\bbL^-$, and that $|\wh\rho(\xi)| <
  1$ for $\xi \neq 0$. Since $\phi$ is analytic in the unit disk of $\C$, it
  follows that $\wh w$ is continuous on $\bbL^-_*$ and analytic
  in the interior of $\bbL^-$. Moreover, $|\wh w(\xi)| \le
  \phi(|\wh\rho(\xi)|) \le \phi(\e^{\Im\xi})$, hence $|\wh w(\xi)| =
  \cO(\e^{\Im\xi})$ as $\Im\xi \to -\infty$. These properties imply (see
  \cite{Sch66}, Ch.~VIII) that $\wh w$ is the Fourier transform of a uniquely
  determined distribution $w \in \cD'(\R)$ with support in $[1,\infty)$.
  
  The injectivity of $\cN$ follows from the facts that the mapping
  $\phi:\set{z}{|z|<1}\to \C$ is locally injective (as $\phi'(z) =
  1/(1{-}Q(z))\neq 0$) and that $\phi:[0,1)\to \R$ is globally injective (as
  $\phi'(s)\geq 1$ for $s\in[0,1)$).  If $\cN(\rho_1)=\cN(\rho_2)$, then, by
  the above, we have $\phi(\wh \rho_1(\xi))=\phi(\wh \rho_2(\xi))$ for $\xi\in
  \bbL^-_*$. This proves $\wh\rho_1(-\i p)= \wh\rho_2(-\i p)$ for $p>0$, as
  $\wh\rho_j(-\i p)\in [0,1)$. By continuity of $\wh\rho_j$ and local
  invertibility we obtain $\wh\rho_1=\wh\rho_2$ on $\bbL^-_*$, and hence
  $\rho_1=\rho_2$.
 
To prove $w=\cN(\rho)\in L^1_\loc([1,\infty))$, choose any $\epsilon > 0$ and
consider the distribution $w_\epsilon:x\mapsto \e^{-\epsilon x}w(x)$. It 
belongs to $\cS'(\R)$ (the space of tempered distributions) and its Fourier 
transform satisfies
\[
  \wh w_\epsilon(\xi) = \wh w(\xi {-}\i\epsilon) = 
  \phi(\wh\rho(\xi {-}\i\epsilon)) \FOR  \Im\xi \le 0.
\]
Now we observe that $\wh\rho(\xi {-} \i\epsilon) = \wh\rho_\epsilon(\xi)$, 
where $\rho_\epsilon(x) = \e^{-\epsilon x}\rho(x)$. Since 
$\|\rho_\epsilon\|_{L^1} \le \e^{-\epsilon} < 1$, the series 
$\sum_{k=1}^\infty \frac{\phi^{(k)}(0)}{k!} \rho_\epsilon^{*k}$ 
converges in $L^1(\R)$ to some function $W_\epsilon \in 
L^1([1,\infty),\R_+)$. (Here we use the crucial fact that 
$\phi^{(k)}(0) \ge 0$ for all $k \in \N$.) By construction, 
\[
  \wh W_\epsilon(\xi) = \sum_{k=1}^\infty \frac{\phi^{(k)}(0)}{k!}
  \bigl(\wh\rho_\epsilon(\xi)\bigr)^k = \phi(\wh\rho(\xi {-} \i\epsilon))
  = \wh w_\epsilon(\xi) \FOR  \Im\xi\le 0, 
\]
giving $w_\epsilon = W_\epsilon \in L^1((1,\infty),\R_+)$, and hence
$w: x\mapsto  \e^{\epsilon x} w_\epsilon(x)$ lies in 
$L^1_\loc([1,\infty),\R_+)$. 
\end{proof}

\begin{remarks}\\
  {\bf 1.} Under the assumptions of Proposition~\ref{S:Ndef}, one has that
  $w=\cN(\rho) \in \cS'(\R)$, i.e., $w$ is a tempered distribution. In fact,
  there exists a constant $C > 0$ such that $|\wh w(\xi)| =
  |\phi(\wh\rho(\xi))| \le C\max\{1,-\log|\xi|\}$ for $\xi \neq 0$, see the
  proof of Proposition~\ref{prop:C.N} below. This means that the singularity
  of $\wh w(\xi)$ at $\xi = 0$ is (not worse than) logarithmic.
  \\[1mm]
  {\bf 2.} More information on $\cN$ can be extracted from the proof of
  Proposition~\ref{S:Ndef}. For instance, if $\rho \in \bbP$, then
  $\cN(\rho)(x) = \rho(x)$  for a.a.\ $x \in (1,n{+}1)$, 
  where 
\begin{equation}\label{eq:ndef} n =
  \min\bigset{j \in \{1,\dots,N\}}{p_j > 0} \ge 1 
\end{equation} 
  is the largest integer such that $|Q(z)| = \cO(|z|^n)$ as $z \to 0$.
  Indeed, in view of \Ref{eq:Uphi}, one has $\phi(z) = z + \cO(|z|^{n+1})$ as
  $z \to 0$.  It follows that
\[
   W_\epsilon = \rho_\epsilon + \sum_{k=n+1}^\infty 
   \frac{\phi^{(k)}(0)}{k!} \,\rho_\epsilon^{*k},  
\]
where the second term in the right-hand side is supported in 
the interval $[n{+}1,\infty)$. Thus $W_\epsilon = \rho_\epsilon$ a.e. 
in $[1,n{+}1]$, which proves the claim. Similarly, using the 
observation that $\supp(\rho_\epsilon^{*k}) \subset [k,\infty)$, 
it is easy to show that, if $\rho : [1,\infty) \to \R_+$ is
continuous, so is $\cN(\rho)$. 
\end{remarks}

The formula \Ref{eq:U.3} is very nice, but does not provide an effective
method for solving the Cauchy problem associated with \Ref{eq:U.1}. 
Indeed, Proposition~\ref{S:Ndef} does not give a sufficient 
characterization of the set $\cN(\bbP)$, which is also the domain of 
$\cN^{-1}$. It is not even clear a priori that this set in left 
invariant by the linear evolution $U(t)$. For this reason, we 
shall use standard PDE techniques to prove existence of solutions 
to \Ref{eq:U.1} in the next section. But the representation 
\Ref{eq:U.3} will be very useful to find self-similar solutions 
of \Ref{eq:U.1} in Section~\ref{SS}, and to study their stability in 
Section~\ref{GC}.


\section{The Cauchy problem for the rescaled system}
\label{s:S}

The evolution equation \Ref{eq:U.1} is not autonomous, and it is 
defined on the time-dependent domain $\set{x \in \R_+ }{ x \ge t }$. 
These drawbacks are eliminated if we rescale the density $\rho(t,x)$ 
by setting
\begin{equation}\label{eq:Srhoeta}
  \rho(t,x) = \frac{1}{t} \eta(\log t, x/t) \FOR  x \ge t \ge 1, 
\end{equation}
or equivalently
\begin{equation}\label{eq:Setarho}
 \eta(\tau,y) = \e^\tau\rho(\e^\tau,\e^\tau y) \FOR  \tau \ge 0, ~y \ge 1. 
\end{equation}
In what follows, we denote by $\tau = \log t$ and $y = x/t$ the new time 
and space coordinates. The rescaled density $\eta(\tau,\cdot)$ now 
belongs to the fixed space $\bbP$ defined in \Ref{eq:Pdef}. 
Moreover, it satisfies the autonomous evolution equation
\begin{equation}\label{eq:S.1}
 \pl_\tau\eta(\tau,y)= \pl_y(y\,\eta(\tau,y)) + \beta(\tau)
 \bbQ[\eta(\tau,\cdot)](y{-}1) \FOR  y \ge 1,
\end{equation}
where $\beta(\tau)=\eta(\tau,1)$ is the new trace which relates to
$\alpha(t)$ via $\beta(\tau)=\e^\tau\alpha(\e^\tau)$. The initial
condition for \Ref{eq:S.1} is $\eta(0,y) = \eta_0(y)$, where $\eta_0 
= \rho_1 \in \bbP$. 

The nonlinearity in \Ref{eq:S.1} has the form $\beta(\tau)  
T_1 \bbQ[\eta(\tau)]$, where $T_1 : \bbP \to \bbP$ is the shift 
operator defined by
\begin{equation}\label{eq:shift}
  (T_1 \eta)(y) = \left\{\begin{array}{cl} \eta(y{-}1) &\text{if }
  y \ge 2;\\0 &\text{if }y < 2.\end{array}\right.
\end{equation}
In particular, for all $\eta \in \bbP$, the support of $T_1
\bbQ[\eta]$ is contained in $[2,\infty)$, or even in $[n{+}1,\infty)$, 
where $n \ge 1$ is defined in \Ref{eq:ndef}. Thus, any solution
of \Ref{eq:S.1} satisfies the linear equation $\pl_\tau\eta = 
\pl_y(y\eta)$ in the strip $\set{(\tau,y)}{ \tau \ge 0,\,1\le y 
\le 2 }$. It follows that $\eta(\tau,y) = \e^{\tau-\tau_0}
\eta(\tau_0,\e^{\tau-\tau_0}y)$ for all $\tau \ge \tau_0 \ge 0$ 
and all $y \ge 1$ such that $\e^{\tau-\tau_0}y \le 2$. Setting 
$y = 1$, we obtain the important relation
\begin{equation}\label{eq:tradef}
  \beta(\tau) = \e^{\tau-\tau_0} \,\eta(\tau_0,\e^{\tau-\tau_0})
  \FOR 0 \le \tau {-} \tau_0 \le \log 2,
\end{equation}
which means that the trace $\beta(\tau)$ for $\tau \in [\tau_0,
\tau_0 {+} \log 2]$ can be determined from the solution 
$\eta(\tau_0,\cdot)$. This formula will be useful to define the 
trace $\beta$ properly when the solution $\eta(\tau,\cdot)$ 
of \Ref{eq:S.1} is not 
continuous. For instance, if $\eta(\tau,\cdot) \in \bbP$ for all 
$\tau \ge 0$ and if $\beta$ satisfies \Ref{eq:tradef}, then 
$\beta \in L^1_\loc([0,\infty),\R_+)$. 

The main purpose of this section is to show that \Ref{eq:S.1} 
defines a well-posed evolution in the space $\bbP$. To do this, 
we consider the associated integral equation
\begin{equation}\label{eq:Sinteq}
  \eta(\tau) = S_\tau \eta_0 + \int_0^\tau \beta(s) S_{\tau-s} 
  T_1 \bbQ[\eta(s)]\dd s \FOR  \tau \ge 0, 
\end{equation}
where $(S_\tau)_{\tau\ge 0}$ is the linear semigroup on $\bbP$ 
defined by 
\begin{equation}\label{eq:S.3}
  (S_\tau \eta)(y) = \left\{\begin{array}{cl}\e^\tau \eta(\e^\tau y) 
  &\text{if } y\ge 1;\\0 &\text{if }y < 1.\end{array}\right.
\end{equation}
To formulate our convergence results in Section~\ref{GC}, we shall
need some weighted $L^p$ spaces which we now introduce. 
For $p \in [1,\infty)$ and $\gamma \ge 0$, we denote by $L^p_\gamma$ 
the function space
\begin{equation}\label{eq:Lpdef}
  L^p_\gamma = \set{w\in L^1_\loc([1,\infty),\R)}{
  \|w\|_{p,\gamma} < \infty },
\end{equation}
where
\[
 \|w\|_{p,\gamma} = \|y^\gamma w\|_{L^p} = \left(\int_1^\infty (y^\gamma
 |w(y)|)^p \dd y\right)^{1/p}.    
\]
When $\gamma = 0$, we simply write $L^p$ instead of $L^p_0$ and 
$\|w\|_p$ instead of $\|w\|_{p,0}$. Remark that $L^p_\gamma 
\hookrightarrow L^1$ if and only if $\gamma > 1-1/p$ (when $p > 1$)
or $\gamma \ge 0$ (when $p = 1$). In what follows, we shall often 
restrict ourselves to such values of $p,\gamma$.

We first give a few basic estimates on the semigroup $(S_\tau)$ and the
nonlinearity $\bbQ$ acting on $L^p_\gamma$.

\begin{lemma}\label{S:semig}
Let $p \in [1,\infty)$ and $\gamma \ge 0$. Then \Ref{eq:S.3} defines
a strongly continuous semigroup $(S_\tau)_{\tau \ge 0}$ in $L^p_\gamma$, 
and
\begin{equation}\label{eq:semig}
   \|S_\tau \eta\|_{p,\gamma} \le \e^{-\tau(\gamma-1+1/p)}
   \|\eta\|_{p,\gamma}, 
\end{equation}
for all $\eta \in L^p_\gamma$ and all $\tau \ge 0$. Moreover, 
equality holds in \Ref{eq:semig} if and only if $\eta(y) = 0$ 
for almost all $y \in [1,\e^\tau]$. 
\end{lemma}

\begin{lemma}\label{S:nonlin}
Let $\bbQ$ be the nonlinear map defined by \Ref{eq:nonlin}.\\
{\bf a)} If $\eta \in L^1$, then $\bbQ[\eta] \in L^1$ and
$\|\bbQ[\eta]\|_1 \le Q(\|\eta\|_1)$. If $\eta,\tilde\eta \in L^1$, 
then 
\[
  \|\bbQ[\eta] - \bbQ[\tilde\eta]\|_1 \le Q'(r) \|\eta -
  \tilde\eta\|_1,
\] 
where $r = \max\{\|\eta\|_1,\|\tilde\eta\|_1\}$. Finally, if $\eta\in
\bbP$, then $\bbQ[\eta] \in \bbP$. \\[1mm]
{\bf b)} Let $p \in [1,\infty)$ and $\gamma > 1-1/p$. If $\eta \in 
L^p_\gamma$, then $\bbQ[\eta] \in L^p_\gamma$, and there exists 
$C > 0$ (independent of $\eta$) such that
\begin{equation}\label{eq:nlin1}
  \|T_1 \bbQ[\eta]\|_{p,\gamma} \le C Q'(\|\eta\|_1) 
  \|\eta\|_{p,\gamma}. 
\end{equation}
If $\eta, \tilde \eta \in L^p_\gamma$ and $R = \max\{\|\eta\|_{p,\gamma},
\|\tilde\eta\|_{p,\gamma}\}$, then
\[
  \|T_1 \bbQ[\eta] - T_1 \bbQ[\tilde \eta]\|_{p,\gamma} \le C Q'(R)
  \|\eta-\tilde\eta\|_{p,\gamma}.
\]
\end{lemma}

\begin{proof}
Estimate \Ref{eq:semig} is a straightforward calculation, 
and the proof of Lemma~\ref{S:nonlin} will be outlined in 
Appendix~\ref{app:C}. 
\end{proof}

We are now ready to state the main result of this section:

\begin{theorem}\label{S:Sexist}
For any $\eta_0 \in L^1((1,\infty),\R)$ with $\|\eta_0\|_1 \le 1$, 
equations \Ref{eq:Sinteq}, \Ref{eq:tradef} have a unique global
solution $\eta \in C^0([0,\infty),L^1)$, which satisfies 
$\|\eta(\tau)\|_1 \le 1$ for all $\tau \ge 0$. In addition, \\
1) if $\eta_0 \in \bbP$, then $\eta(\tau) \in \bbP$ for all 
$\tau \ge 0$;\\
2) if $\eta_0 \in L^p_\gamma$ for some $p \ge 1$ and some $\gamma >
1-1/p$, then $\eta \in C^0([0,\infty),L^p_\gamma)$. 
\end{theorem}

\begin{proof} Fix $\eta_0 \in B_1$, where $B_1 = \set{\eta \in L^1 }{
\|\eta\|_1 \le 1 }$. Setting $\tau_0 = 0$ in \Ref{eq:tradef}, we
obtain
\begin{equation}\label{eq:traceid}
  \beta(\tau) = \e^\tau \eta_0 (\e^\tau) \FOR  0 \le \tau \le \log 2.
\end{equation}
The first step is to show that \Ref{eq:Sinteq},\Ref{eq:traceid} 
have a unique solution $\eta \in C^0([0,\log 2],L^1)$. 

Let $q = Q'(1) \ge 1$, and let $T = (\log 2)/m$, where $m \in \bbN^*$
is sufficiently large so that, for all $k = 1, \dots,m$, 
\begin{equation}\label{eq:mcond}
  \int_{(k-1)T}^{kT} \e^s |\eta_0(\e^s)| \dd s <
  \frac{1}{q}. 
\end{equation}
We introduce the Banach space $X = C^0([0,T],L^1)$ equipped with 
the norm
\[
   \|\eta\|_X = \sup_{0 \le \tau \le T}\|\eta(\tau)\|_1. 
\]
Let $B = \set{\eta \in X}{\|\eta\|_X \le 1 }$, and let $F : X 
\mapsto X$ be the nonlinear map defined by
\[
   (F[\eta])(\tau) =  S_\tau \eta_0 + \int_0^\tau \beta(s) S_{\tau-s} 
  T_1 \bbQ[\eta(s)]\dd s \FOR  0 \le \tau \le T, 
\]
where $\beta(s)$ is given by \Ref{eq:traceid}. We claim that $F(B) 
\subset B$ and that $F$ is a strict contraction in $B$. Indeed:

\noindent {\bf a)} Assume that $\eta \in B$. Using Lemmas~\ref{S:semig} 
and \ref{S:nonlin}, we find, for all $\tau \in [0,T]$, 
\begin{eqnarray}\nonumber
  \|(F[\eta])(\tau)\|_1 &\le & \|S_\tau \eta_0\|_1 + \int_0^\tau
    |\beta(s)| \|S_{\tau-s}T_1 \bbQ[\eta(s)]\|_1 \dd s \\ \nonumber
  &=& \int_1^\infty \e^\tau |\eta_0(\e^\tau y)|\dd y +  \int_0^\tau
    |\beta(s)| \| T_1 \bbQ[\eta(s)]\|_1 \dd s \\ \label{eq:Farray}
  &=& \int_{\e^\tau}^\infty |\eta_0(y)|\dd y + \int_0^\tau
    \e^s |\eta_0(\e^s)| \|\bbQ[\eta(s)]\|_1 \dd s \\ \nonumber
  & \le & \int_{\e^\tau}^\infty |\eta_0(y)|\dd y + Q(\|\eta\|_X)
    \int_1^{\e^\tau} |\eta_0(y)|\dd y \ \le\ 1, 
\end{eqnarray}
since $Q(\|\eta\|_X) \le Q(1) = 1$ and $\|\eta_0\|_1 \le 1$. 
This shows that $F(B) \subset B$. 

\noindent {\bf b)} If $\eta, \tilde \eta\in B$, then for all 
$\tau \in [0,T]$,
\begin{eqnarray*}
  \|(F[\eta])(\tau) {-} (F[\tilde\eta])(\tau)\|_1 &\le & \int_0^\tau
    |\beta(s)| \|S_{\tau-s}(T_1 \bbQ[\eta(s)] {-} T_1 
    \bbQ[\tilde\eta(s)])\|_1 \dd s \\ 
  & = & \int_0^\tau \e^s |\eta_0(\e^s)| \|\bbQ[\eta(s)] {-} 
    \bbQ[\tilde\eta(s)])\|_1 \dd s \\
  &\le & \int_0^\tau \e^s |\eta_0(\e^s)| Q'(1) \|\eta(s) {-}
    \tilde\eta(s)\|_1 \dd s \\
  &\le& q \Bigl(\int_0^T \e^s |\eta_0(\e^s)| \dd s \Bigr) 
    \|\eta{-}\tilde\eta\|_X. 
\end{eqnarray*}
In view of \Ref{eq:mcond}, this shows that $F$ is a strict contraction
in $B$. 

Let $\eta \in X$ be the unique fixed point of $F$ in the ball $B$.
Then $\eta$ satisfies \Ref{eq:Sinteq}, and using Gronwall's lemma it
is readily verified that $\eta$ is in fact the unique solution of 
\Ref{eq:Sinteq} in the whole space $X = C^0([0,T],L^1)$. 
Repeating the same argument $m$ times (where $m$ is such that 
\Ref{eq:mcond} holds), we conclude that equations~\Ref{eq:Sinteq},
\Ref{eq:traceid} have a unique solution $\eta \in
C^0([0,\log 2],L^1)$, which satisfies $\|\eta(\tau)\|_1 \le 1$ 
for all $\tau \in [0,\log 2]$. Moreover, it is clear that 
\Ref{eq:tradef} holds for all $\tau_0 \in [0,\log 2]$ and almost all 
$\tau \in [\tau_0,\log 2]$. 

For $\tau \in [0,\log 2]$, let $\Xi_\tau : B_1 \to B_1$ be the 
nonlinear map defined by $\Xi_\tau \eta_0 = \eta(\tau)$, where 
$\eta(\tau)$ is the solution of \Ref{eq:Sinteq} we have just constructed. 
Then it is easy to verify that
$\Xi_{\tau_1+\tau_2} = \Xi_{\tau_1} \circ \Xi_{\tau_2}$ for 
$0 \le \tau_1 + \tau_2 \le \log 2$. It follows that the family 
$(\Xi_\tau)$ can be extended to a continuous semiflow 
$(\Xi_\tau)_{\tau\ge 0}$. By construction, if $\eta_0 \in B_1$ and
if we set $\eta(\tau) = \Xi_\tau \eta_0$ for all $\tau \ge 0$, 
then $\eta \in C^0([0,\infty),L^1)$ is the unique solution 
of \Ref{eq:Sinteq}, \Ref{eq:tradef}, and $\eta(\tau) \in B_1$
for all $\tau \ge 0$. This proves the first part of Theorem~\ref{S:Sexist}. 

Assume now that $\eta_0 \in \bbP$. Keeping the same notations as
above, we define 
\[
  \tilde B = \set{\eta \in X }{ \eta(\tau) \in \bbP \text{ for all }
  \tau \in [0,T] }.
\]
In particular, $\tilde B$ is a closed subset of $B$, as $\bbP$ is closed in
$B_1\subset L^1$. If $\eta \in 
\tilde B$, it is clear that $(F[\eta])(\tau) \in L^1((1,\infty),\bbR_+)$ 
for all $\tau \in [0,T]$, and that all inequalities in \Ref{eq:Farray}
can be replaced by equalities. Thus $F(\tilde B) \subset \tilde B$, 
hence the solution $\eta \in C^0([0,\infty),L^1)$ of
\Ref{eq:Sinteq} satisfies $\eta(\tau) \in \bbP$ for all $\tau \in 
[0,T]$. Proceeding as above, we then show that $\eta(\tau) \in \bbP$ 
for all $\tau \in [0,\log 2]$, hence for all $\tau \ge 0$. This
proves assertion 1) in Theorem~\ref{S:Sexist}. 

Finally, assume that $\eta_0 \in L^p_\gamma$ for some $p \ge 1$ and
some $\gamma > 1-1/p$, and that $\|\eta_0\|_1 \le 1$. Using 
Lemmas~\ref{S:semig}, \ref{S:nonlin} and a fixed point argument as
before, it is straightforward to show that the solution 
$\eta \in C^0([0,\infty),L^1)$ of \Ref{eq:Sinteq} satisfies 
$\eta \in C^0([0,T],L^p_\gamma)$ for some $T > 0$ (depending on 
$\eta_0$). Let
$$
   T^* = \sup\bigset{T > 0}{\eta \in C^0([0,T],L^p_\gamma)} 
   \in (0,\infty].
$$
We claim that $T^* = \infty$. Indeed, assume on the contrary 
that $0 < T^* < \infty$. Since $\|\eta(\tau)\|_1 \le 1$ for all 
$\tau \ge 0$, it follows from \Ref{eq:Sinteq}, \Ref{eq:semig}, 
\Ref{eq:nlin1} that 
\[
   \|\eta(\tau)\|_{p,\gamma} \le \|\eta_0\|_{p,\gamma} +
   C q \int_0^\tau |\beta(s)| \|\eta(s)\|_{p,\gamma} \dd s \FOR  
   0 \le \tau < T^*. 
\]
Using Gronwall's lemma and the fact that $\beta \in L^1_\loc([0,\infty))$, 
we deduce that $\|\eta(\tau)\|_{p,\gamma} \le C'$ for all $\tau \in 
[0,T^*)$. In view of \Ref{eq:Sinteq}, \Ref{eq:nlin1}, this 
in turn implies that $\eta(\tau)$ has a limit in $L^p_\gamma$ 
as $\tau \nearrow T^*$, giving $\eta \in C^0([0,T^*],L^p_\gamma)$.
Since we have a local existence result in $L^p_\gamma$, 
we conclude that $\eta \in C^0([0,T],L^p_\gamma)$ for some 
$T > T^*$, which contradicts the definition of $T^*$. This proves 
assertion 2) in Theorem~\ref{S:Sexist}.
\end{proof}

The nonlinear map $\cN$ introduced in the previous section can also be 
used to linearize \Ref{eq:S.1}. Indeed, the Fourier transforms of 
$\rho$ and $\eta$ are related via $\wh\rho(t,\xi) = \wh\eta(\log t, t\xi)$,
so that \Ref{eq:Srhoeta} is just a rescaling of the Fourier variable
$\xi$. As is clear from \Ref{eq:Ndef}, this transformation commutes 
with the action of $\cN$. Thus, if $\rho$ is a solution of
\Ref{eq:U.1} with initial data $\rho_1$ and if $\eta$ is the
corresponding solution of \Ref{eq:S.1} given by \Ref{eq:Setarho}, 
it follows from \Ref{eq:U.3} that
\begin{equation}\label{eq:undo}
  \frac1t \cN\bigl(\eta(\log t,\cdot)\bigr)(x/t) = 
  \cN(\eta_0)(x) \FOR  x \ge t \ge 1, 
\end{equation}
where $\eta_0 = \rho_1$. Setting $\tau = \log t$ and $y = x/t$, 
we obtain the representation formula
\begin{equation}\label{eq:explicit}
  \cN(\eta(\tau)) = S_\tau\, \cN(\eta_0) \FOR  \tau \ge 0,
\end{equation}
where $(S_\tau)$ is the linear semigroup \Ref{eq:S.3}.
The last result of this section shows that this formula is indeed
correct:

\begin{proposition}\label{S:resol}
Let $\eta_0 \in \bbP$, and let $\eta \in C^0([0,\infty),\bbP)$ 
be the solution of \Ref{eq:Sinteq} given by Theorem~\ref{S:Sexist}. 
Then $\cN(\eta(\tau)) = S_\tau \,\cN(\eta_0)$ for all $\tau \ge 0$. 
\end{proposition}

\begin{proof} We establish the formula by returning to the unscaled variables
  $(t,x)$ and by showing that the formal steps of Section \ref{unscaled} 
  can be made rigorous for the solutions of \Ref{eq:S.1}. 
Define $\rho : [1,\infty)^2 \to \R_+$ by $\rho(t,x) = \frac1t 
\eta(\log t,x/t)$ if $x \ge t \ge 1$ and $\rho(t,x) = 0$ if 
$1 \le x < t$. Then $\rho \in C^0([1,\infty),\bbP)$, and 
rescaling \Ref{eq:Sinteq} we find
\begin{equation}\label{eq:newint}
  \rho(t) = U(t)\left(\rho_1 + \int_1^t \alpha(s) T_s\bbQ[\rho(s)]
  \dd s \right) \FOR  t \ge 1,
\end{equation}
where $\rho_1 = \eta_0 \in \bbP$, $\alpha(t) = \frac1t \beta(\log t)$, 
$U(t)$ is the linear operator \Ref{eq:Udef}, and $T_s$ is the shift
operator defined as in \Ref{eq:shift}. To simplify the notation, 
we set $f(s,x) = (T_s\bbQ[\rho(s)])(x)$. Then $f \in C^0([1,\infty),
L^1)$, so that $(s,x) \mapsto \alpha(s)f(s,x) \in L^1_\loc([1,\infty),
L^1)$. By construction, the trace $\alpha$ satisfies the identity
\[
  \alpha(t) = \rho_1(t) + \int_1^t \alpha(s) f(s,t)\dd s \FOR  
  \ae t \ge 1. 
\]

We now apply the Fourier transform to \Ref{eq:newint}. For any 
$\xi \in \bbL^-$ and any $t \ge 1$, we find 
\[
  \wh\rho(t,\xi) = \int_t^\infty \rho_1(x) \e^{-\i \xi x} \dd x
  + \int_t^\infty \Big\{\int_1^t \alpha(s)f(s,x)\dd s \Big\}
  \,\e^{-\i \xi x} \dd x.
\]
Since $\rho_1 \in \bbP$, the first term in the right-hand side is
absolutely continuous with respect to $t$, and
\[
  \pl_t \int_t^\infty \rho_1(x) \e^{-\i \xi x} \dd x = 
  -\rho_1(t) \e^{-\i \xi t} \FOR \ae t \ge 1. 
\]
The second term can be decomposed as $h_1(t,\xi)-h_2(t,\xi)$, 
where
\begin{eqnarray*}
  h_1(t,\xi) &=& \int_1^\infty \Big\{\int_1^t \alpha(s)f(s,x)\dd s 
   \Big\} \,\e^{-\i \xi x} \dd x \,=\, \int_1^t \alpha(s) \wh f(s,\xi)
   \dd s, \\
  h_2(t,\xi) &=& \int_1^t \Big\{\int_1^t \alpha(s)f(s,x)\dd s 
   \Big\} \,\e^{-\i \xi x} \dd x.  
\end{eqnarray*}
Clearly, $h_1(t,\xi)$ is absolutely continuous with respect to $t$, 
and
\[
  \pl_t h_1(t,\xi) = \alpha(t) \wh f(t,\xi) = \alpha(t) 
   \e^{-\i \xi t} Q(\wh\rho(t,\xi)) \FOR \ae t \ge 1. 
\]
Next, since $f(s,x) = 0$ for $x < s$, we have $\int_1^t
\alpha(s)f(s,x)\dd s = \int_1^x \alpha(s)f(s,x)\dd s$, 
and this expression is a locally integrable function of $x$. 
It follows that $h_2(t,\xi)$ is absolutely continuous with respect 
to $t$, and
\[
  \pl_t h_2(t,\xi) = \e^{-\i \xi t} \int_1^t \alpha(s)f(s,t)\dd s 
  \FOR \ae t \ge 1. 
\]
Summarizing, we have shown that, for any $\xi \in \bbL^-$, the
Fourier transform $\wh\rho(t,\xi)$ is absolutely continuous with 
respect to $t$ and satisfies
\begin{eqnarray*}
  \pl_t \wh\rho(t,\xi) &=& -\e^{-\i \xi t} \left(\rho_1(t) + 
   \int_1^t \alpha(s)f(s,t) \dd s\right) + \alpha(t) \e^{-\i \xi t}
  Q(\wh\rho(t,\xi)) \\
   &=& \alpha(t)\,\e^{-\i \xi t}\Big(Q(\wh\rho(t,\xi))-1\Big) \quad
  \FOR \ae t \ge 1.
\end{eqnarray*}
This gives \Ref{eq:U.2}. Now, proceeding exactly as in Section~\ref{unscaled}, 
we deduce that \Ref{eq:Wsol} holds for all $t \ge 1$ if $\Im\xi < 0$, 
and this in turn is equivalent to \Ref{eq:U.3}. Finally, using the
transformation \Ref{eq:undo} we obtain \Ref{eq:explicit}. 
\end{proof}


\section{Properties of the steady states}\label{SS}

This section is devoted to the time-independent solutions of
\Ref{eq:S.1} in the space $\bbP$ defined by \Ref{eq:Pdef}. 

\begin{definition}
We say that $\eta_0 \in \bbP$ is a {\it steady state} of \Ref{eq:S.1}
if the solution $\eta \in C^0([0,\infty),\bbP)$ of \Ref{eq:Sinteq}
given by Theorem~\ref{S:Sexist} satisfies $\eta(\tau) = \eta_0$ for 
all $\tau \ge 0$. 
\end{definition}

\noindent
The steady states of \Ref{eq:S.1} will also be called ``equilibria'' 
or ``stationary solutions''. 

\begin{lemma}\label{S:beta}
If $\eta_0 \in \bbP$ is a steady state of \Ref{eq:S.1}, there
exists $\beta \geq 0$ such that $\eta_0(y) = \beta/y$ for almost all
$y \in [1,2]$. 
\end{lemma}

\begin{proof}
If $\eta(\tau) \equiv \eta_0$, \Ref{eq:Sinteq} implies that 
$\eta_0(y) = \e^\tau \eta_0(\e^\tau y)$ for all $\tau \in [0,\log 2]$
and  $\ae y \in [1,2\,\e^{-\tau}]$, because the nonlinearity
in \Ref{eq:Sinteq} vanishes identically for such values of $\tau,y$.
We define $F:x\mapsto \int_1^{\e^x}\eta_0(y)\dd y \geq 0$ and obtain
\[
  F(x{+}y)=F(x)+F(y)\FOR x,y\geq 0 \AND x{+}y\leq \log 2.
\]
Since $F$ is continuous, we conclude that $F(x)=\beta x$ for some
$\beta \geq 0$. 
Differentiating implies $\beta=\e^x \eta_0(\e^x)$ for $\ae x\in
[0,\log 2]$ which gives the desired result. 
\end{proof}

Let $\eta_0 \in \bbP$ be a steady state. Since $\eta_0$ coincides
almost everywhere in $[1,2]$ with a continuous function, the constant 
$\beta$ in Lemma~\ref{S:beta} can be identified with  $\eta_0(1)$. 
Clearly, the trace function defined by \Ref{eq:tradef} satisfies 
$\beta(\tau) = \beta$ for all $\tau \ge 0$. In particular, the 
integral equation \Ref{eq:Sinteq} reduces to
\begin{equation}\label{eq:intstat}
  \eta_0 = S_\tau \eta_0 + \beta \int_0^\tau S_s T_1\bbQ[\eta_0]\dd s 
  \FOR \tau \ge 0.
\end{equation}
From $\eta_0\in \bbP$ we now conclude that $\beta>0$.

On the other hand, if $\eta_0 \in \bbP$ and $w = \cN(\eta_0)$, it 
follows from Propositions~\ref{S:Ndef} and
\ref{S:resol} that $\eta_0$ is a steady state if and only if 
$S_\tau w = w$ for all $\tau \ge 0$. In view of \Ref{eq:S.3}, this is 
the case if and only if there exists $\beta' \in \R$ such that 
$w = \beta' w^*$, where 
\begin{equation}\label{eq:wstar}
  w^*(y) = \left\{\begin{array}{cl} 1/y &\text{if }
  y \ge 1,\\0 &\text{if }y < 1.\end{array}\right.
\end{equation}
But since $w(y) = \eta_0(y)$ for $\ae y \in [1,2]$ (see Remark 2 after Proposition~\ref{S:Ndef}), we necessarily have
$\beta' = \beta = \eta_0(1)$. 

Finally, since equilibria are time-independent solutions 
of \Ref{eq:S.1}, we certainly expect them to solve the 
ordinary differential equation
\begin{equation}\label{eq:time-indep}
  (y\eta)'(y) + \beta (T_1 \bbQ[\eta])(y) = 0 \FOR  y \ge 1,
  \quad \eta(1) = \beta. 
\end{equation}
Remark that the initial value $\beta$ also appears as a parameter 
in front of the nonlinear term. It is not difficult to show that
\Ref{eq:time-indep} has global solutions:

\begin{lemma}\label{Seta}
For any $\beta \in \R$, equation \Ref{eq:time-indep} has a unique
global solution $\eta : [1,\infty) \to \R$. 
\end{lemma}

\begin{proof}
For any $k \in \N$, let $I_k = [kn+1,(k{+}1)n+1]$, where $n \in \N_*$ 
is defined in \Ref{eq:ndef}. For any $\eta \in L^1_\loc([1,\infty),\R)$, 
the nonlinear term $(T_1 \bbQ[\eta])(y)$ only depends on the values of 
$\eta(z)$ for $z \le y-n$. In particular, $(T_1 \bbQ[\eta])(y) = 0$ for 
$y \le n+1$, so that any solution of \Ref{eq:time-indep} satisfies
$\eta(y) = \beta/y$ for $y \in I_0 = [1,n+1]$. Using this 
information, one can compute $(T_1 \bbQ[\eta])(y)$ explicitly for
$y \in I_1 = [n+1,2n+1]$, and then solve \Ref{eq:time-indep} on this 
interval to determine $\eta(y)$ for $y \in I_1$. By construction, 
$\eta$ is smooth on both $I_0$ and $I_1$, but $\eta$ has a 
discontinuity of order $n$ at $y = n+1$, in the sense that the 
derivatives $\eta^{(k)}(y)$ are continuous for $k = 0, \dots,n-1$, 
whereas $\eta^{(n)}(y)$ has different limits to the left and to the 
right at $y = n+1$ (if $\beta \neq 0$). Iterating this procedure, we 
find that \Ref{eq:time-indep} has a unique global solution $\eta \in 
C^{n-1,1}([1,\infty),\R)$, which satisfies $\eta \in C^\infty(I_k)$ 
for all $k \in \N$. 
\end{proof}

The following result shows that equilibria of \Ref{eq:S.1}
indeed correspond to solutions of the differential equation
\Ref{eq:time-indep}. 

\begin{proposition}\label{S:equiv}
If $\eta_0 \in \bbP$ and $\beta > 0$, the following assertions are
equivalent:\\
{\bf a)} $\eta_0$ is a steady state of \Ref{eq:S.1} with 
$\eta_0(1) = \beta$.\\
{\bf b)} $\eta_0$ coincides almost everywhere with the solution
of \Ref{eq:time-indep}.\\
{\bf c)} $\cN(\eta_0) = \beta w^*$. 
\end{proposition}

\begin{proof}
We already proved that ${\bf a)} \Leftrightarrow {\bf c)}$. 
If $\eta_0 \in \bbP$ is a steady state with $\eta_0(1) = \beta$, 
it follows from \Ref{eq:intstat} that
\[
   \frac{S_\tau \eta_0 -\eta_0}{\tau} \,+\, \frac{\beta}{\tau} 
   \int_0^\tau S_s T_1\bbQ[\eta_0]\dd s \,=\, 0,  
\]
for all $\tau > 0$. Using \Ref{eq:S.3}, it is not difficult
to verify that the first term converges to $(y\eta_0)'$ in 
$\cD'((1,\infty))$ as $\tau \to 0$, while the second one tends
to $\beta T_1 \bbQ[\eta_0]$ in $L^1((1,\infty))$. This shows
that (after modification on a set of measure zero) $\eta_0$ is 
absolutely continuous on $(1,\infty)$ and satisfies the differential 
equation \Ref{eq:time-indep} for almost all $y > 1$. It follows 
easily that $\eta_0$ is the solution of \Ref{eq:time-indep} in the 
sense of Lemma~\ref{Seta}. Thus ${\bf a)} \Rightarrow {\bf b)}$.

Conversely, assume that $\eta_0 \in \bbP$ satisfies \Ref{eq:time-indep}.
Applying the semi-group $S_\tau$ to \Ref{eq:time-indep} and
integrating over $\tau$, we immediately obtain \Ref{eq:intstat}, 
which implies that $\eta_0$ is a steady state. This proves that
${\bf b)} \Rightarrow {\bf a)}$.
\end{proof}

The main goal of this section is to determine for which values 
of $\beta > 0$ the solution $\eta$ of \Ref{eq:time-indep} 
belongs to $\bbP$. Our strategy is to use the characterization
${\bf c)}$ in Proposition~\ref{S:equiv}. Therefore, we are led to 
study the image of $\beta w^*$ under the map $\cN^{-1}$, and this 
requires very precise information on the complex transformations 
\Ref{eq:Uphi} and \Ref{eq:Upsi}. The following quantities, related 
to the polynomial $Q(z)$, will play an important role in the sequel:
\begin{equation}\label{eq:kappadef}
  q = Q'(1) \ge 1\ \AND \   
  \kappa = \exp\left( \int_0^1 \Bigl( \frac{1}{1-z} -\frac{q}{1-Q(z)}
  \Bigr) \dd z\right) \le 1. 
\end{equation}

\begin{lemma}\label{Sphi}
Let
\[
   \Phi(z) = 1 - \e^{-q\phi(z)} \FOR  |z| < 1, 
\]
where $\phi$ is defined in \Ref{eq:Uphi}. Then $\Phi$ can be extended
analytically to a neighborhood of the real positive axis $\R_+$. This 
extension satisfies $\Phi(z) \ge 0$ and $\Phi'(z) > 0$ for all $z \ge 0$. 
Moreover, $\Phi(0) = 0$, $\Phi'(0) = q$, $\Phi(1) = 1$, $\Phi'(1) = 
\kappa$, and $\Phi(z) \to R$ as $z \to \infty$, where
\begin{equation}\label{eq:Rdef}
   R = 1 + \exp\left( \int_0^2 \Bigl( \frac{1}{1-z} -\frac{q}{1-Q(z)}
  \Bigr)\dd z - \int_2^\infty \frac{q}{1-Q(z)} \dd z\right). 
\end{equation}
Note that $R = \infty$ if $Q(z) = z$ and $1 < R < \infty$ otherwise. 
\end{lemma}

\begin{proof}
Since the polynomial $1 - Q(z)$ has the unique real positive root $z = 1$, 
which is a simple root because $Q'(1) = q \neq 0$, it is clear that 
the function
\[
   \chi(z) = \exp\left( \int_0^z \Bigl( \frac{1}{1-t} -\frac{q}{1-Q(t)}
  \Bigr) \dd t\right) = \frac{\e^{-q\phi(z)}}{1-z} \FOR  |z| < 1, 
\]
can be extended to an analytic map in a neighborhood of the real 
positive axis $\R_+$. Moreover, $\chi(0) = 1$, $\chi(1) = \kappa$, 
and using $z{-}1=\exp(-\int_2^z\frac{\d t}{1-t})$ 
shows that $(z{-}1)\chi(z) \to R{-}1$ for $z \to 
\infty$, where $R$ is defined in \Ref{eq:Rdef}. Since $\Phi(z) = 
1 - (1{-}z)\chi(z)$, we conclude that the function $\Phi$ has the desired 
properties. In particular, 
\[
   \Phi'(z) = q \chi(z) \frac{1-z}{1-Q(z)}, 
\]
so that $\Phi'(z) > 0$ for all $z \ge 0$. 
\end{proof}

It follows from Lemma~\ref{Sphi} that the map $\Phi : [0,\infty) \to 
[0,R)$ is one-to-one and onto. Let $\Psi = \Phi^{-1} : [0,R) \to 
[0,\infty)$ be the inverse map. Then $\Psi(0) = 0$, $\Psi'(0) = 1/q$, 
$\Psi(1) = 1$, $\Psi'(1) = 1/\kappa$, and $\Psi'(u) > 0$ for all 
$u \in [0,R)$. By construction, 
\begin{equation}\label{eq:U2psi}
   \Psi(u) = \psi\Bigl(-\frac{1}{q} \log(1{-}u)\Bigr) \FOR 
   0 \le u < 1. 
\end{equation}

\begin{lemma}\label{Spsi}
The function $\Psi : [0,R) \to [0,\infty)$ is absolutely monotone, 
i.e. $\Psi^{(k)}(u) \ge 0$ for all $k \in \N$ and all $u \in [0,R)$.
In particular, $\Psi$ can be extended to an analytic function 
on the disc $|u| < R$, and there exist nonnegative coefficients 
$(\Psi_k)_{k\in\N_*}$ such that
\[
   \Psi(u) = \sum_{k=1}^\infty \Psi_k u^k \FOR  |u| < R. 
\]
\end{lemma}

\begin{proof}
Since $\Psi = \Phi^{-1}$, we already know that $\Psi$ is analytic in 
a neighborhood of $[0,R)$. We first show by induction that, 
for all $n \in \N_*$, there exists a polynomial $P_n$ such that
\begin{equation}\label{eq:Urec1}
  \Psi^{(n)}(u) = \frac{P_n(\Psi(u))}{q^n \,(1{-}u)^n}\quad\text{for }
   0 < u < 1.
\end{equation}
Indeed, differentiating \Ref{eq:U2psi} and using \Ref{eq:Upsi}, we obtain
\begin{equation}\label{eq:Psidiff}
  \Psi'(u) = \frac{1-Q(\Psi(u))}{q\,(1{-}u)} \quad\text{for } 0 < u < 1. 
\end{equation}
Thus \Ref{eq:Urec1} holds for $n = 1$ with $P_1(z) = 1-Q(z)$. On the 
other hand, differentiating \Ref{eq:Urec1} and using \Ref{eq:Psidiff}, 
we find, for $0<u<1$,
\begin{equation}\label{eq:Prec}
  \Psi^{(n+1)}(u) = \frac{P_{n+1}(\Psi(u))}{q^{n+1} \,(1{-}u)^{n+1}} 
  \quad \text{with } P_{n+1}(z) = P_n'(z)(1{-}Q(z)) + nq P_n(z). 
\end{equation}
Therefore, \Ref{eq:Urec1} is established.

We next show that, for all $n \in \N_*$, there exists a polynomial $R_n(z)$
with {\sl nonnegative coefficients} such that
\begin{equation}\label{eq:Urec2}
  P_n(z) = (1{-}Q(z)) (1{-}z)^{n-1} R_n(z). 
\end{equation}
Obviously, \Ref{eq:Urec2} holds for $n = 1$ with $R_1(z) = 1$. 
Combining \Ref{eq:Prec} and \Ref{eq:Urec2}, we obtain the recursion 
relation
\[
   R_{n+1}(z) = A_1(z) R_n'(z) + A_2(z) R_n(z) + (n{-}1)A_3(z) R_n(z), 
\]
where the coefficient functions  $A_j$ are given by
\begin{eqnarray*}
  A_1(z) &=& \frac{1-Q(z)}{1-z} = \sum_{j=1}^N p_j \frac{1-z^j}{1-z}, \\
  A_2(z) &=& \frac{q-Q'(z)}{1-z} = \sum_{j=2}^N j p_j \frac{1-z^{j-1}}{1-z},\\
  A_3(z) &=& \frac{q}{1-z} - \frac{1-Q(z)}{(1-z)^2} = 
    \sum_{j=2}^N p_j \sum_{k=1}^{j-1} \frac{1-z^k}{1-z}.
\end{eqnarray*}
Because of $p_j\geq 0$ all $A_1, A_2, A_3$ are polynomials (in $z$) with
nonnegative coefficients. Thus, the same property holds for $R_n$ by induction
over $n$.

Since $0 < \Psi(u) < 1$ and $0 < Q(\Psi(u)) < 1$ for all $ u \in(0,1)$, 
it follows from \Ref{eq:Urec1} and \Ref{eq:Urec2} that $\Psi^{(n)}(u) \ge 0$
for all $n \in \N$ and all $u \in (0,1)$, hence also for $u \in [0,1]$. 
By a classical result of Bernstein (see \cite{Fe71}, Section VII.2), 
the power series
\begin{equation}\label{eq:Useries}
  \sum_{k=1}^\infty \Psi_k u^k, \quad \text{where } \Psi_k = \frac{1}{k!}
  \Psi^{(k)}(0) \ge 0,
\end{equation}
converges absolutely and uniformly for $|u| \le 1$, and defines an 
analytic continuation of $\Psi$ to the unit disk. Moreover, if $R_1 
\ge 1$ denotes the radius of convergence of the series
\Ref{eq:Useries}, it is well-known (see for instance \cite{Ru87},
exercise 16.1) that the analytic function defined by \Ref{eq:Useries} 
has a singularity at $u = R_1$. Since $\Psi(u) \to \infty$ as $u \nearrow
R$, it follows that $R = R_1$. This concludes the proof. 
\end{proof}

\begin{example}
To conclude this study of the mappings $\Phi$ and $\Psi$, we give 
an explicit example of a nonlinearity $Q$ for which these
functions can be calculated explicitly. Let $Q(z)=(1{-}a)z+az^2$,
where $a\in[0,1]$. The value $a=1$ corresponds to the coarsening
equation \Ref{eq:roeq}, while $a=0$ is a particular case of a model
studied in \cite{CaP00}. Then $q = 1+a = 1/\kappa$, $R = 1 + 1/a$, 
and
\[
  \phi(z)=\frac{1}{1+a}\log\frac{1+az}{1-z},\qquad
  \psi(w)=\frac{1-\e^{-qw}}{1+a\e^{-qw}}.
\]
The auxiliary functions $\Phi$, $\Psi$ are:
\[
 \Phi(z) = \frac{(1{+}a)z}{1+az}, \qquad
 \Psi(u) = \frac{u}{1{+}a-au}.
\]
\end{example}

We are now ready to state and prove the main result of this section.

\begin{theorem}\label{thm:equi} {\rm (Steady states of \Ref{eq:S.1})}\\
Fix $\theta > 0$ and let $\eta^*_\theta : [1,\infty) \to \R$ be 
the solution of \Ref{eq:time-indep} with $\beta = \theta/q$.
Then\\
{\bf a)} $\eta^*_\theta \in \bbP$ if and only if $0 < \theta \le 1$.\\
{\bf b)} If $\theta \in (0,1]$, $\eta^*_\theta \in \bbP$ is positive 
and strictly decreasing, so that $y\eta^*_\theta(y) \to 0$ as 
$y \to \infty$. \\
{\bf c)} If $0 < \theta < 1$, then 
\begin{equation}\label{eq:asym1} 
  \lim_{y \to\infty} y^{1+\theta} \eta^*_\theta(y) = \frac{\theta \,\e^{\theta
  \euler}}{\kappa \,\Gamma(1{-}\theta)},
\end{equation}
where $\Gamma$ is the Gamma function and $\euler = -\Gamma'(1) 
\approx 0.577216$ is Euler's constant.\\
{\bf d)} If $\theta = 1$, then
\begin{equation}\label{eq:mom1} 
  \int_1^\infty y\eta^*_1(y)\dd y = \frac{\e^{\euler}}{\kappa}.
\end{equation}
Moreover, if\/ $\deg Q > 1$, there exists $\lambda > 0$ such that
\begin{equation}\label{eq:asym2} 
  \lim_{y \to\infty} \frac{\log \eta^*_1(y)}{y} = -\lambda.
\end{equation}
For $Q(z) = z$ we have
\begin{equation}\label{eq:asym3} 
  \lim_{y \to\infty} \frac{\log \eta^*_1(y)}{y\log y} = -1.
\end{equation}
\end{theorem}

\begin{remark}
It follows from Theorem~\ref{thm:equi} and Proposition~\ref{S:equiv}
that \Ref{eq:S.1} has a {\it unique} steady state $\eta^*_1 \in \bbP$ 
such that $\int_1^\infty y \eta^*_1(y)\dd y < \infty$.
\end{remark}

\begin{proof}
We first show that $\eta^*_\theta \in \bbP$ if $0 < \theta \le 1$. 
According to Proposition~\ref{S:equiv}, it is sufficient to prove
that there exists an element of $\bbP$ (still denoted by 
$\eta^*_\theta$) such that $\cN(\eta^*_\theta) = (\theta/q)w^*$. 
Since $\cN^{-1} =  \cF^{-1} \circ \psi \circ \cF$ and $\psi(w) = 
\Psi(1{-}\e^{-qw})$ by \Ref{eq:U2psi}, this relation is equivalent to
\begin{equation}\label{eq:etathdef}
  \wh \eta^*_\theta = \Psi\bigl(1{-} \e^{-\theta \wh w^*} \bigr), 
\end{equation}
where $\wh \eta^*_\theta = \cF \eta^*_\theta$ and $\wh w^* = \cF w^*$. 
In view of \Ref{eq:wstar},
\begin{equation}\label{eq:E2exp}
  \wh w^*(\xi) = \int_1^\infty \frac{\e^{-\i \xi y}}{y}\dd y = 
  \Eone(\i \xi), 
\end{equation}
where $\Eone$ is the exponential integral, see \cite{AS72}. It is
well-known that
\begin{equation}\label{eq:E1exp}
  \Eone(z) = -\log z - \euler + \chi(z) \FOR  \left|\arg z\right| < \pi, 
\end{equation}
where $\chi : \bbC \to \bbC$ is an entire function with $\chi(0) =
0$ and $\chi'(0) = 1$. Thus, $\wh w^*$ is analytic in the interior 
of $\bbL^-$, where $\bbL^- = \set{\xi \in \bbC }{\Im \xi \le 0}$. 
Moreover, $\Re(\wh w^*(\xi)) \to \infty$ as $\xi \to 0$ within $\bbL^-$. 

In Appendix~\ref{app:A}, we prove that $|1-\e^{-\theta \wh w^*(\xi)}| 
 < 1$ for all $\xi \in \bbL^- \setminus \{0\}$ and all $\theta\in (0,1]$, 
see also Figure~\ref{figure2}. From 
Lemma~\ref{Spsi}, we also know that $\Psi$ is analytic in the disk 
of radius $R > 1$ centered at the origin. Therefore, the map 
$\wh \eta^*_\theta$ defined by \Ref{eq:etathdef} is continuous over
$\bbL^-$ (with $\wh \eta^*_\theta(0) = 1$) and analytic in the
interior of $\bbL^-$. In addition, since $|\Psi(u)| \le |u|$ whenever 
$|u| \le 1$, we have the bound
\[
   |\wh \eta^*_\theta(\xi)| \le |1{-}\e^{-\theta \wh w^*(\xi)}| 
   \le 2\theta |\wh w^*(\xi)| \FOR \xi \in \bbL^- \setminus \{0\}. 
\]
In particular, $|\wh \eta^*_\theta(\xi)| = \cO(\e^{\Im \xi})$ as
$\Im \xi \to -\infty$. By the Paley-Wiener Theorem (see for instance 
\cite{Ru87}), we conclude that $\eta^*_\theta = \cF^{-1} \wh \eta^*_\theta
\in L^2((1,\infty))$.
 
To prove that $\eta^*_\theta$ is nonnegative, we argue as in \cite{CaP92}.
Consider the Laplace transform $\wt \eta^*_\theta=\cL \eta^*_\theta$, 
which satisfies $\wt \eta^*_\theta(p) = \wh \eta^*_\theta(-\i p)$. As is
well-known (see \cite{Fe71}, Section XIII.4), positivity of
$\eta^*_\theta$ is equivalent to {\it complete monotonicity} of 
$\wt\eta^*_\theta$, namely $(-1)^k {\wt\eta^*_\theta}{}^{(k)}(p) \ge 0$ 
for all $k\in \N$ and $p > 0$. Recall that 
\begin{equation}\label{eq:wteta} 
  \wt\eta^*_\theta(p) = \Psi(1{-}\e^{-\theta \wt w^*(p)}) = 
  \Psi(1{-} \e^{-\theta \Eone(p)}) \FOR  p > 0. 
\end{equation}
We apply Lemma~\ref{lem:comp} below with 
\[
  f_1:\left\{\ba{ccc} (0,1) &\to& \R, \\ u &\mapsto& \Psi(1{-}u);\ea\right.
  \quad \AND \quad g_1:\left\{
  \ba{ccc} (0,\infty) &\to& (0,1), \\ p &\mapsto&  
  \e^{-\theta \Eone(p)}. \ea \right.
\] 
By Lemma~\ref{Spsi}, $f_1$ is completely monotone, thus it remains to 
show that $g_1'$ is completely monotone. Observe that $g_1' = f_2 
\circ g_2$, where $f_2 : \R \to \R$ is defined by $f_2(w) = 
\theta \e^{-w}$ and $g_2 : (0,\infty) \to \R$ 
by
\[
  g_2(p) = \theta \Eone(p) - \log (-\Eone'(p)) = \theta \Eone(p) + p +\log p.
\]
Clearly, $f_2$ is completely monotone, thus (again by Lemma~\ref{lem:comp}) 
it remains to prove that $g'_2$ is completely monotone. This follows 
from the representation
\[
  g'_2(p) = -\theta \frac{\e^{-p}}p + 1 +\frac1p = 1 + (1{-}\theta)\frac1p
  +\theta \int_0^1 \e^{-sp}\dd s.
\]
Thus, we have shown that $\eta^*_\theta \in L^2((1,\infty),\R_+)$. Since 
$\wt \eta(p) \to 1$ as $p \searrow 0$, we conclude that $\eta^*_\theta 
\in L^1$ and $\int_1^\infty \eta^*_\theta(y)\dd y = 1$, i.e.,
$\eta^*_\theta \in \bbP$. 

Now, fix $\theta > 1$ and assume that $\eta^*_\theta \in \bbP$, where
$\eta^*_\theta$ is the solution of \Ref{eq:time-indep} with $\beta =
\theta/q$. According to Proposition~\ref{S:equiv}, $\cN(\eta^*_\theta)
= (\theta/q)w^*$, so that \Ref{eq:wteta} holds. Thus, in view of
\Ref{eq:E1exp}, the Laplace transform of $\eta^*_\theta$  
satisfies
\begin{equation}\label{eq:Lapeta} 
  \wt\eta^*_\theta(p) = \Psi\Bigl(1 - p^\theta\,\e^{\theta(\euler-\chi(p))}
  \Bigr) = 1 - \kappa^{-1}p^\theta\,\e^{\theta \euler} + \cO(p^{1+\theta})
  \FOR p \searrow 0. 
\end{equation}
Since $\theta > 1$, it follows that $\int_1^\infty y\eta^*_\theta(y)\dd y=
-(\wt\eta^*_\theta)'(0) =0$, which clearly contradicts the hypothesis
$\eta^*_\theta \in \bbP$. This proves {\bf a)}.

Next, fix $\theta \in (0,1]$. To prove that $\eta^*_\theta$ is strictly 
decreasing, it is sufficient to show that $\eta^*_\theta(y) > 0$ 
for all $y \ge 1$, since $y(\eta^*_\theta)'(y) + \eta^*_\theta(y) \le 0$
by \Ref{eq:time-indep}. Assume on the contrary
that there exists $y_0 > 1$ such that $\eta^*_\theta(y_0) = 0$ and 
$\eta^*_\theta(y) > 0$ for $1 \le y < y_0$. It is clear that $y_0 >
n{+}1$, where $n$ is defined in \Ref{eq:ndef}. Thus, 
$(T_1 \bbQ[\eta])(y_0) > 0$, hence $\eta^*_\theta{}'(y_0) < 0$ by 
\Ref{eq:time-indep}, which contradicts the fact that $\eta^*_\theta 
\in \bbP$. This proves {\bf b)}. 

Assume now that $0 < \theta < 1$. In Appendix~\ref{app:B}, we prove
that the limit in the left-hand side of \Ref{eq:asym1} exists. 
Let $L(\theta)$ denote this limit, and let
\[
   H_\theta(y) = \int_y^\infty \eta^*_\theta(x) \dd x \FOR  y \ge 1.
\] 
Clearly, $y^\theta H_\theta(y) \to L(\theta)/\theta$ as $y \to
\infty$. Thus, the Laplace transform of $H_\theta$ satisfies
\[
  p^{1-\theta} \wt H_\theta(p) = \int_p^\infty \e^{-t}t^{-\theta} \,
  \Bigl(\frac{t}{p}\Bigr)^\theta H_\theta\Bigl(\frac{t}{p}\Bigr)\dd t
  \to \Gamma(1{-}\theta)\frac{L(\theta)}{\theta} \quad 
  \text{as } p \searrow 0. 
\]
Since $\wt \eta^*_\theta(p) = \e^{-p} - p\wt H_\theta(p) = 1 - p^\theta 
\Gamma(1{-}\theta)L(\theta)/\theta + \mbox{\scriptsize $\cO$}
(p^\theta)$ as $p \searrow 0$, 
it follows from \Ref{eq:Lapeta} that $\Gamma(1{-}\theta)L(\theta)/\theta
= \e^{\theta\euler}/\kappa$. This proves \Ref{eq:asym1}. 

Finally, let $\theta = 1$. Then \Ref{eq:wteta}, \Ref{eq:Lapeta} show 
that the Laplace transform $\wh \eta^*_1$ is analytic in the half-plane 
$\set{p \in \C }{\Re p > -\lambda }$, where $\lambda > 0$ is the
unique real root of the equation $1 - \e^{-\Eone(-\lambda)} = R$ (if 
$Q(z) = z$, then $R = \infty$, hence also $\lambda = \infty$.) 
In particular, $\eta^*_1(y)$ decays exponentially as $y \to \infty$, and 
\[
  -{\wt\eta^*_1}{}'(0) = \int_1^\infty y\eta^*_1(y)\dd y = 
  \frac{\e^{\euler}}{\kappa}. 
\]
If $\deg Q > 1$, then $\lambda < \infty$, and the arguments given 
in \cite{CaP92} (in the particular case $Q(z) = z^2$) show that 
\Ref{eq:asym2} holds. If $Q(z) = z$, then $\lambda = \infty$ and 
$\eta^*_1(y) = \rho(y-1)/y$, where $\rho : [0,\infty) \to \R_+$ is 
the Dickmann function studied in \cite{CaP00}. From the asymptotics
of $\rho$ given there, we deduce that \Ref{eq:asym3} holds. 
This concludes the proof.
\end{proof}

\figurewithtex{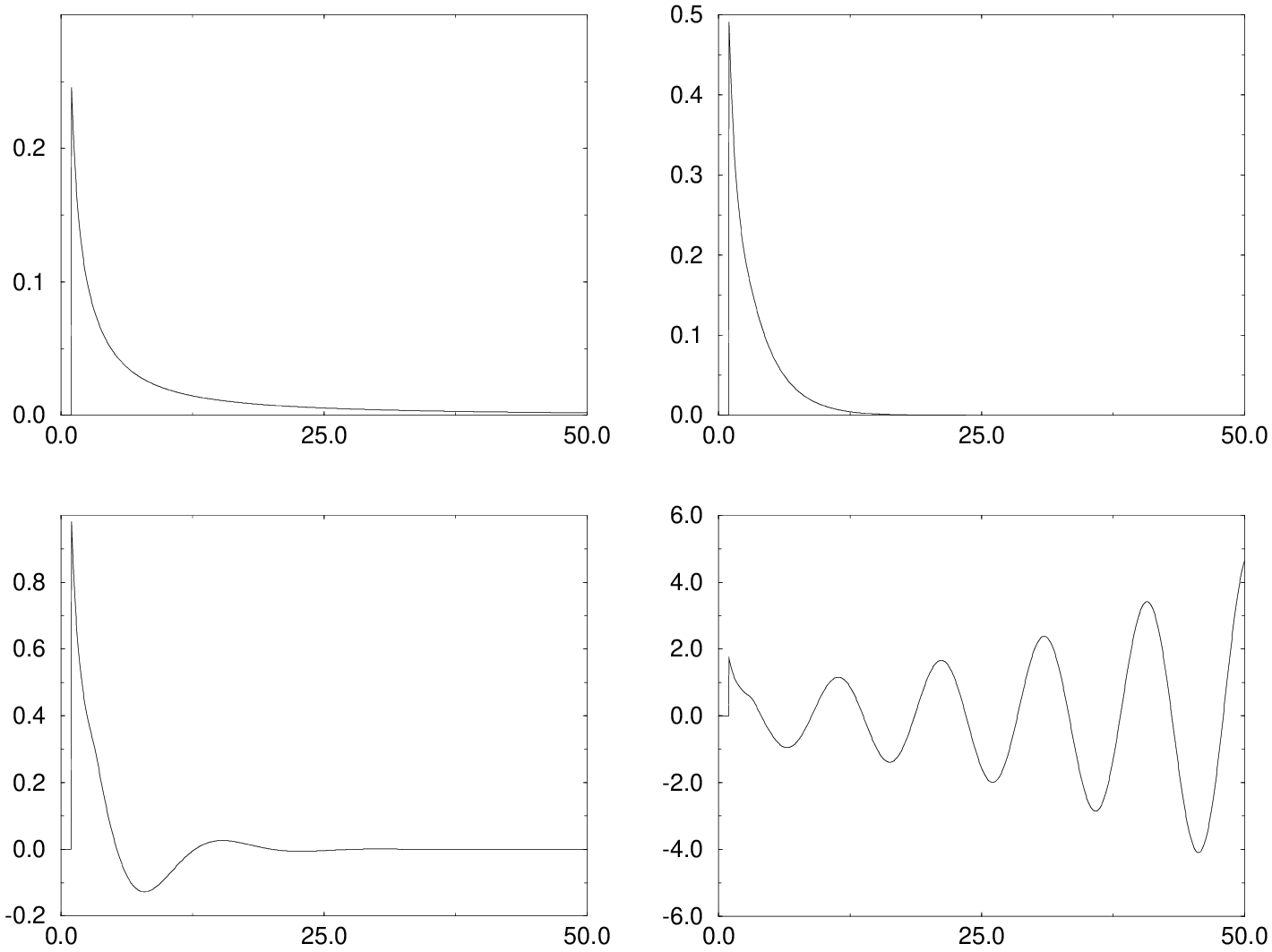}{fig1.tex}{11}{14}{figure1}{%
The steady state  $\eta_\theta^*$ of the coarsening equation \Ref{eq:S.1} with
$Q(z) = z^2$ is represented for four values of the parameter
$\theta$. The first two graphs ($\theta = 0.5$ and $\theta = 1.0$) 
illustrate the conclusions of Theorem~\ref{thm:equi}, and the other 
two ($\theta =2.0$ and $\theta = 3.5$) the remarks after 
Lemma~\ref{lem:comp}. The pictures were produced using  
the explicit formula \Ref{eq:sssol} and a FFT routine to compute
the Fourier transforms.} 

The following lemma was used in the proof of Theorem~\ref{thm:equi}. For its
proof see  \cite{Fe71}, Section XIII.4.  

\begin{definition}
Let $I \subset \R$ be an open interval, 
and let $f \in C^\infty(I,\R)$. The function $f$ is called {\it completely
monotone} if $(-1)^k f^{(k)}(x) \ge 0$ for all $x \in I$ and 
all $k \in \N$.  
\end{definition}

\begin{lemma}\label{lem:comp} 
{\rm (Composition lemma)}\\
Let $I,J \subset \R$ be open intervals. If $f : J \to \R$ is
completely monotone and $g:I \to J$ has a derivative $g'$ which is 
completely monotone, then $f \circ g : I \to \R$ is completely monotone.
\end{lemma}

\begin{remarks}
The (generalized) steady states $\eta^*_\theta$ with $\theta > 1$ will 
not be studied in this paper, because they do not lie in our function 
space $\bbP$. We just mention here a few properties that can 
established using the techniques developed in the proof of 
Theorem~\ref{thm:equi}. There exists a critical value 
$\theta_* \in (1,\infty]$ such 
that\\[1mm]
{\bf 1)} If $1 < \theta < \theta_*$, then $\eta^*_\theta \in 
L^1((1,\infty),\R)$ and $\int_1^\infty \eta^*_\theta(y)\dd y = 1$. 
However, $\eta^*_\theta$ is not a positive function. In particular, 
$\|\eta^*_\theta\|_1 > 1$, so that $\eta^*_\theta$ does not belong to the 
unit ball of $L^1$ where existence of global solutions is known from 
Theorem~\ref{S:Sexist}.\\[1mm]
{\bf 2)} If $\theta > \theta_*$, then $\eta^*_\theta \notin 
L^1((1,\infty),\R)$.\\[1mm]
Moreover, $\theta_* = \infty$ if $Q(z) = z$, whereas $\theta_* 
< \infty$ if $\deg Q > 1$. In the particular case where 
$Q(z) = z^2$, one has $\theta_* \approx 3.24826$. These statements 
are illustrated in Figure~\ref{figure1}.  
\end{remarks}


\section{Global convergence results}
\label{GC}

In this final section, we use the explicit representation formula 
\Ref{eq:explicit} to 
study the long-time behavior of the solutions of \Ref{eq:S.1}. In
particular, we obtain global stability results for the steady states 
$\eta^*_\theta$ with $0 < \theta \le 1$. 

Since the nonlinear map $\cN$, which allows us to linearize 
\Ref{eq:S.1}, has a simple expression in Fourier variables, it is
convenient to use $L^2$--based function spaces instead of the $L^1$--based
function spaces which are more natural for the existence theory. 
Our basic space will be 
\[
   P_\gamma = \bbP \cap L^2_\gamma \FOR  \gamma \ge 0,
\]
where $\bbP$ is defined in \Ref{eq:Pdef} and $L^2_\gamma$ in 
\Ref{eq:Lpdef}. Remark that $P_\gamma$ is a closed subspace 
of $L^2_\gamma$ if $\gamma > 1/2$, since $L^2_\gamma \hookrightarrow
L^1$. 

The image of $P_\gamma$ under the Fourier transform $\cF$ can be 
characterized completely. Let $\bbH^\gamma_1$ be the space of 
all functions $z : \bbL^-\to \C$ satisfying the following three 
conditions: \\[1mm]
(i) $z$ is analytic in the interior of $\bbL^-$,\\[1mm]
(ii) for each $\xi_2 \le 0$, the map $\xi_1 \mapsto z(\xi_1{+}\i\xi_2)$ 
lies in the Sobolev space $H^\gamma(\R)$,\\[1mm]
(iii) 
\begin{equation}\label{eq:bbHnorm}
  \|z\|_{\bbH^\gamma_1} = \sup_{\xi_2\leq 0}\e^{-\xi_2}\|z(\cdot{+}\i 
  \xi_2)\|_{H^\gamma(\R)} < \infty.
\end{equation}
(In \Ref{eq:bbHnorm}, the supremum over $\xi_2 \le 0$ is always
attained at $\xi_2=0$.) 

Then $\eta \in L^2_\gamma$ if and only if $\wh\eta = \cF \eta
\in \bbH^\gamma_1$ (when $\gamma = 0$, this is just the Paley-Wiener
theorem, see \cite{Ru87}; the general case follows using the Fourier
characterization of the Sobolev space $H^\gamma(\R)$.) Moreover, 
the map $\eta \mapsto \|\wh \eta\|_{\bbH^\gamma_1} = 
\|\wh \eta\|_{H^\gamma(\R)}$ is a norm on $L^2_\gamma$ which is 
equivalent to $\|\eta\|_{2,\gamma}$. If in addition $\eta \in
P_\gamma$, then $\wh \eta(0) = 1$ and $\xi_1 \mapsto \wh \eta(\xi_1)$ 
is a positive definite function on $\R$ (in the sense of Bochner). 
Furthermore, 
$|\wh \eta(\xi)| < 1$ for all $\xi \in \bbL^- \setminus\{0\}$. 

Assume now that $\eta \in P_\gamma$ for some $\gamma > 1/2$, 
and let $w = \cN(\eta)$, namely $\wh w(\xi) = \phi(\wh \eta(\xi))$. 
Since $\phi$ is analytic in the unit disk, it is clear that 
$\wh w$ is analytic in the interior of $\bbL^-$. Moreover, the 
fact that $\phi(z)= z + \cO(|z|^2)$ as $z \to 0$ guarantees that 
$\wh w(\xi)$ has the same decay properties as $\wh \eta(\xi)$ as 
$|\xi| \to \infty$. However, since $\wh \eta(0) = 1$ and since
$\phi(z)$ has a singularity at $z = 1$, we see that $\wh w(\xi)$ 
necessarily has a singularity at $\xi = 0$. This is the reason 
why the nonlinear transformation $\cN$ {\it does not} map 
$P_\gamma$ into itself. To handle this difficulty, 
our strategy is to subtract from $\wh w(\xi)$ a suitable function
with the same singularity at $\xi = 0$ and whose inverse
Fourier transform is explicitly known. 

If $\gamma > 3/2$, a natural candidate for this counter-term is
$\frac{1}{q} \wh w^*(\xi) = \phi(\wh \eta^*_1(\xi))$, where $\eta^*_1$ is 
the unique steady state of \Ref{eq:S.1} that belongs to $L^2_\gamma$, 
see Theorem~\ref{thm:equi}. We recall that $w^*$ is defined in 
\Ref{eq:wstar}. 

\begin{proposition}\label{prop:C.N}
Let $\gamma > 3/2$ and $\eta\in P_\gamma$. Then $\cN(\eta) =
\frac{1}{q} w^* + d$ with $d \in L^2_{\gamma-1}$.
\end{proposition}

\begin{proof}
We first show that $\cN(\eta) \in L^2 = L^2((1,\infty))$. 
As explained above, it is sufficient to prove that $\wh w \equiv
\phi(\wh \eta)$ satisfies \Ref{eq:bbHnorm} with $\gamma = 0$. 
Choose $a > 0$ sufficiently small so that $|\phi(z)| \le 2|z|$ for all 
$z \in \C$ with $|z| \le a$. Since $\eta \in L^2_\gamma$, 
there exists $b < 0$ such that $|\wh \eta(\xi)| \le a$ whenever
$\Im\xi \le b$. Thus
\[
  \sup_{\xi_2\le b}\e^{-\xi_2}\|\wh w(\cdot{+}\i \xi_2)\|_{L^2}
  \le 2 \sup_{\xi_2\le b}\e^{-\xi_2}\|\wh \eta(\cdot{+}\i
  \xi_2)\|_{L^2} < \infty.
\]
On the other hand, by a variant of the Riemann-Lebesgue lemma, 
there exists $c > 0$ such that $|\wh \eta(\xi)| \le a$ for all 
$\xi \in \bbL^-$ with $|\Re \xi| \ge c$. Arguing as before, we 
thus get
\[
  \sup_{b \le \xi_2\le 0} \int_{|\xi_1| \ge c} 
  |\wh w(\xi_1{+}\i \xi_2)|^2 \dd \xi_1 < \infty~.
\]
It remains to verify that 
\begin{equation}\label{eq:whbound}
  \int_{|\xi_1| \le c} |\wh w(\xi_1{+}\i \xi_2)|^2 \dd \xi_1 \le C 
  \quad \text{uniformly in } \xi_2 \in [b,0].
\end{equation}
Since $\wh w : \bbL^- \to \C$ is continuous except at the origin, 
it is sufficient to establish \Ref{eq:whbound} for $b,c$ sufficiently
small. Now, as $\xi \to 0$ in $\bbL^-$, we have the expansion
\[
  \wh\eta(\xi) = 1-\i \mu \xi + r_1(\xi), \quad \text{with }
  r_1(\xi) = \left\{
  \begin{array}{lcl} \cO(|\xi|^2) & {\rm if} & \gamma > 5/2,\\ 
  \cO(|\xi|^{\gamma-1/2}) & {\rm if} & 3/2 < \gamma <  5/2,
  \end{array}\right.
\]
where $\mu = \int_1^\infty y\eta(y)\dd y > 1$. Using the
representation $\phi(z) = -(1/q)\log(1{-}\Phi(z))$ together with 
the properties of $\Phi$ listed in Lemma~\ref{Sphi}, we thus obtain
\[
  \wh w(\xi) = -\frac{1}{q} \log(1{-}\Phi(\wh \eta(\xi))) = 
  -\frac{1}{q} \log(\i\kappa\mu\xi {+} r_2(\xi)), 
\]
where $\kappa = \Phi'(1)$ and $r_2(\xi)$ satisfies the same bounds 
as $r_1(\xi)$. This expansion immediately implies \Ref{eq:whbound}
if $b,c$ are sufficiently small. Thus, we have shown that $\cN(\eta) 
\in L^2$. Since obviously $\frac{1}{q} w^* \in L^2$, we deduce that 
$d = \cN(\eta) - \frac{1}{q} w^* \in L^2$, too. 

To prove that $d \in L^2_{\gamma-1}$, it remains to verify that 
$\wh d \in H^{\gamma-1}(\R)$. Again, by a localization argument, 
it is sufficient to show that $\wh d \in H^{\gamma-1}((-c,c))$ 
for some $c > 0$ sufficiently small. If $\xi \in \R$, $|\xi| < c$, 
we use the representation
\begin{equation}\label{eq:repr}
  \wh d(\xi) = \phi(\wh \eta(\xi)) -  \frac{1}{q}\wh  w^*(\xi)
  = -\frac{1}{q} \log\Bigl(\frac{1-\Phi(\wh \eta(\xi))}
     {\e^{-\wh w^*(\xi)}}\Bigr).
\end{equation}
From \Ref{eq:E1exp}, we know that $\e^{-\wh w^*(\xi)} = \i\xi
\e^{\euler}\e^{-\chi(\i\xi)}$, where $\chi$ is an entire 
function vanishing at the origin. It follows that $\wh d(\xi) = 
-(1/q)\log(D(\xi)/\xi)$, where $D \in H^\gamma((-c,c))$ satisfies
$D(0) = 0$ and $D'(0) = \kappa\mu\e^{-\euler}$. The claim is 
now a direct consequence of Lemma~\ref{lem:Hardy} below. 
This concludes the proof of Proposition~\ref{prop:C.N}.
\end{proof}

\begin{remark}
It follows immediately from the proof of Proposition~\ref{prop:C.N} 
that
\begin{equation}\label{eq:momid1}
  \wh d(0) = \int_1^\infty d(y)\dd y = \frac{1}{q}(\euler - 
  \log(\kappa\mu)), \quad \text{where } \mu = 
  \int_1^\infty y\eta(y)\dd y.  
\end{equation}  
\end{remark}

\begin{lemma}\label{lem:Hardy} Let $\gamma \ge 1$, and let 
$I \subset \R$ be an open interval containing $0$. There exists a 
constant $C(I,\gamma) > 0$ such that, for each $f \in H^\gamma(I)$ 
with  $f(0) = 0$, there exists $g \in H^{\gamma-1}(I)$ such that 
$f(x) = xg(x)$ for all $x \in I$ and 
\[
  \|g\|_{ H^{\gamma-1}(I)} \le C(I,\gamma) \|f\|_{H^\gamma(I)}.
\]
\end{lemma}

\begin{proof} 
It is sufficient to prove the claim for $I = \R$ (the general
case can be reduced to this one using a bounded extension operator).
If $f\in H^\gamma(\R)$ and $f(0)=0$, the Fourier transform 
$\wh f$ has zero mean and satisfies $\lambda^\gamma \wh f 
\in L^2(\R)$, where $\lambda(\xi) = (1{+}\xi^2)^{1/2}$.
Define $g \in L^2(\R)$ by its Fourier transform
\[
   \i \wh g(\xi) = \int_{-\infty}^\xi \wh f(s)\dd s = -\int_\xi^\infty
   \wh f(s)\dd s \FOR  \xi \in \R.
\]
Then $x g(x) = f(x)$ for (almost) all $x \in \R$. Moreover, since
\[
   \lambda(\xi)^{\gamma-1} |\wh g(\xi)| \le \left\{\ba{cc} 
   \int_\xi^\infty \lambda(s)^{\gamma-1} |\wh f(s)|\dd s
     &\text{if }\xi \ge 0,\\[1mm] 
   \int_{-\infty}^\xi \lambda(s)^{\gamma-1} |\wh f(s)|\dd s 
   &\text{if }\xi\le 0,\ea\right.
\]
it follows from Theorem~328 in \cite{HLP59} that $\|\lambda^{\gamma-1}
\wh g\|_{L^2} \le 2\|\lambda^\gamma\wh f\|_{L^2}$, which is the
desired bound.
\end{proof}

We next show that the inverse map $\cN^{-1}$ is well-defined 
in a neighborhood of $\frac{1}{q} w^*$ in $L^2_\gamma$. 

\begin{proposition}\label{prop:C.NN}
Let $\gamma > 1/2$. There exists $\eps > 0$ such that, for all 
$d \in L^2_\gamma$ with $\|d\|_{2,\gamma} \le \eps$, the function
$\cN^{-1}(\frac{1}{q} w^* {+} d)$ is well-defined and lies in $L^2_\gamma$. 
Moreover, there exists $C > 0$ such that
\[
  \|\cN^{-1}({\ts\frac{1}{q}} w^* {+} d) - \eta^*_1\|_{2,\gamma} \le  
  C\|d\|_{2,\gamma},
\]
where $\eta^*_1 = \cN^{-1}(\frac{1}{q} w^*)$. 
\end{proposition}

\begin{proof}
Throughout the proof, we denote by $\|\cdot\|_\gamma$ instead of 
$\|\cdot\|_{\bbH^\gamma_1}$ the norm on $\bbH^\gamma_1$ defined 
by \Ref{eq:bbHnorm}. We first remark that the space $\bbH^\gamma_1$
is an algebra if $\gamma > 1/2$: there exists $C_1 > 0$ such that 
$\|r s\|_\gamma \le C_1 \|r\|_\gamma \|s\|_\gamma$ for all 
$r,s \in \bbH^\gamma_1$. Moreover, as is easy to verify, there 
exists $C_2 > 0$ such that, for all integer $k \ge 1$, 
\begin{equation}\label{eq:power}
  \|r^k\|_\gamma \le C_2 k^{\ell(k,\gamma)} \|r\|_\infty^{k-\ell(k,\gamma)} 
  \|r\|_\gamma^{\ell(k,\gamma)}, 
\end{equation}
where $\|r\|_\infty = \sup\set{|r(\xi)|}{\xi \in \bbL^-} \le
C\|r\|_\gamma$ and $\ell(k,\gamma) = \min\set{n \in \N}{n \ge k \text{ or }
n \ge \gamma}$. 

Assume that $d \in L^2_\gamma$ for some $\gamma > 1/2$, so 
that $\wh d \in \bbH^\gamma_1$. For all $\xi \in \bbL^-$, we 
define
\begin{equation}\label{eq:rsdef}
  r(\xi) = 1 - \e^{-\wh w^*(\xi)} \AND  s(\xi) = \e^{-\wh w^*(\xi)}
  (1 {-} \e^{-q\wh d(\xi)}).  
\end{equation}
From \Ref{eq:E2exp}, \Ref{eq:E1exp}, it is easy 
to see that $r \in \bbH^\gamma_1$, and we prove in Appendix~\ref{app:A}
that $\|r\|_\infty \le 1$. On the other hand, since $\bbH^\gamma_1$
is an algebra, it is clear that $\sigma = 1 {-} \e^{-q\wh d} \in 
\bbH^\gamma_1$, hence $s = (1{-}r)\sigma \in \bbH^\gamma_1$. In addition, 
if $\|d\|_{2,\gamma} \le \eps$ for some $\eps \le 1$, there
exists $C_3 > 0$ such that $\|s\|_\gamma \le C_3 \|d\|_{2,\gamma} 
\le C_3 \eps$. In particular, $\|s\|_\infty \le C\eps$. 

We now fix $R_1 \in (1,R)$, where $R > 1$ is defined in \Ref{eq:Rdef}, 
and we assume that $\eps \le 1$ is sufficiently small so that 
$\|s\|_\infty \le R_1-1$. We then define $\wh \eta \in \bbH^\gamma_1$
by 
\[
  \wh \eta = \psi({\ts\frac{1}{q}}\wh  w^* {+} \wh d) = \Psi(1 {-} 
  \e^{-q(\frac{1}{q}\wh  w^* + \wh d)})
  = \Psi(r{+}s), 
\]
where $\Psi$ is given by \Ref{eq:U2psi}. From Lemma~\ref{Spsi}, we 
know that $\Psi$ is analytic in the disk $\set{u \in \C }{|u| < R }$, 
with the expansion $\Psi(u) = \sum_{k\ge 1} \Psi_k u^k$. 
Since $r+s \in \bbH^\gamma_1$ and $\|r+s\|_\infty \le R_1 < R$, 
it follows from \Ref{eq:power} that the series $\Psi(r{+}s)$ 
converges in $\bbH^\gamma_1$, so that $\wh \eta \in \bbH^\gamma_1$. 
By construction, $\wh \eta = \cF \eta$ for some $\eta \in L^2_\gamma$
with $\cN(\eta) = \frac{1}{q} w^* + d$. 

It remains to show that $\|\wh \eta - \wh \eta^*_1\|_\gamma \le C
\|d\|_{2,\gamma}$, where $\wh \eta^*_1 = \Psi(r)$, see \Ref{eq:wteta}. 
Using 
\[
\|(r{+}s)^{k} - r^{k}\|_\gamma\leq C_1 \sup\bigset{ \|k(r{+}\theta
  s)^{k-1}\|_\gamma }{\theta\in [0,1]}\, \|s\|_\gamma
\]
and \Ref{eq:power}, it is straightforward to verify that there
exists $C_4 > 0$ such that, for all $k \ge 1$, 
\begin{equation}\label{eq:rsbdd}
  \|(r{+}s)^k - r^k\|_\gamma \le C_4 k^{2+\gamma} R_1^{k-1}\|s\|_\gamma.
\end{equation} 
Since
\[
  \wh \eta - \wh \eta^*_1 = \Psi(r{+}s) - \Psi(r) = \sum_{k=1}^\infty 
  \Psi_k ((r{+}s)^k - r^k), 
\]
it follows that
\[
  \|\wh \eta - \wh \eta^*_1\|_\gamma \le C_4 \Bigl(\sum_{k=1}^\infty
  k^{2+\gamma} \Psi_k R_1^{k-1}\Bigr)\|s\|_\gamma \le 
  C_5 \|d\|_{2,\gamma}. 
\]
This concludes the proof. 
\end{proof}

\begin{remark}
Unlike $\cN$, the inverse mapping $\cN^{-1}$ is not positivity 
preserving. However, if in Proposition~\ref{prop:C.NN} we 
{\it assume} in addition that $\eta = \cN^{-1}(\frac{1}{q} w^*{+}d)$ is 
a positive function, then $y \mapsto y\eta(y) \in L^1((1,\infty))$ and
\begin{equation}\label{eq:momid2}
   \int_1^\infty y\eta(y)\dd y = \frac1\kappa \,\e^{\euler-qd_0}, 
   \quad \text{where } d_0 = \int_1^\infty d(y) \dd y.
\end{equation}
Indeed, on the one hand the Laplace transform $\wt \eta(p) = \wh \eta(-\i p)$ 
satisfies
\[
  \frac{1-\wt\eta(p)}p = \frac1p \Bigl(1 - \Psi(1{-}p
  \,\e^{\euler-\chi(p)-q\wt d(p)})\Bigr)
  \  \longrightarrow\ 
  \frac1\kappa \,\e^{\euler-qd_0} \FOR p\searrow 0,
\]
and on the other hand, using $\eta(y) \ge 0$, we find
\[
  \frac{1-\wt\eta(p)}p = \int_1^\infty y\eta(y)
  \frac{1-\e^{-py}}{py}\dd y
 \ \longrightarrow  \
   \int_1^\infty y\eta(y) \dd y=\| \eta \|_{1,1} \FOR p\searrow 0,  
\]
by Lebesgue's monotone convergence theorem.
\end{remark}
 
In addition to Lemma~\ref{S:nonlin}, the following bounds on the 
nonlinearity $\bbQ[\eta]$ will be used to prove our convergence
results:

\begin{lemma}\label{lem:bdd}
Fix $\gamma > 3/2$. For any $M > 0$, there exists $C > 0$ such that
the following estimates hold:\\[1mm]
{\bf a)} For all $\eta \in P_\gamma$ with 
$\|\eta\|_{2,\gamma-1} \le M$ and $\|\eta\|_{1,1} \le M$, 
\begin{equation}\label{eq:bdd1}
  \|T_1\bbQ[\eta]\|_{2,\gamma} \le q\|\eta\|_{2,\gamma} + C.
\end{equation}
{\bf b)} If $\gamma \ge 2$, then for all $\eta,\tilde\eta \in P_\gamma$ 
with $\|\eta\|_{2,\gamma} \le M$ and $\|\tilde\eta\|_{2,\gamma}
\le M$, 
\begin{equation}\label{eq:bdd2}
 \|T_1\bbQ[\eta] {-} T_1\bbQ[\tilde\eta]\|_{2,\gamma} \le
  q\|\eta {-} \tilde\eta\|_{2,\gamma} + C\|\eta {-}
  \tilde\eta\|_{2,\gamma-1}.
\end{equation}
\end{lemma}

\begin{proof}
See Appendix~\ref{app:C}.  
\end{proof}

We are now ready to state the main result of this section, which 
shows that all solutions of \Ref{eq:S.1} in $P_\gamma$ with 
$\gamma > 3/2$ converge towards the limiting profile $\eta^*_1$. 

\begin{theorem}\label{thm:conv}
Assume that $\eta_0 \in P_\gamma$ for some $\gamma>3/2$, and let 
$\eta \in C^0([0,\infty),P_\gamma)$ be the solution of \Ref{eq:Sinteq} 
given by Theorem~\ref{S:Sexist}. Then $\eta$ is bounded in $L^2_\gamma$
and there exists $C > 0$ such that 
\begin{equation}\label{eq:prelim}
  \|\eta(\tau)-\eta^*_1\|_{2,\gamma-1} \le C\,\e^{-(\gamma-3/2)\tau}
  \quad\text{for } \tau \ge 0. 
\end{equation}
Moreover, if $\gamma \ge 2$, then
\[ 
  \|\eta(\tau)-\eta^*_1\|_{2,\gamma} \le C(1{+}\tau)
  \,\e^{-(\gamma-3/2)\tau} \quad\text{for } \tau \ge 0. 
\]
\end{theorem}

Remarkably, the faster the initial data decay at infinity, the
faster the solution converges to the steady state. For compactly 
supported data, it should be possible to obtain faster decay 
than exponential.

\begin{proof}
We first prove \Ref{eq:prelim} using the representation formula 
\Ref{eq:explicit}. By Proposition~\Ref{prop:C.N}, $\cN(\eta_0) = 
\frac{1}{q} w^* + d$ for some $d\in L^2_{\gamma-1}$. Since the
semigroup $S_\tau$ is linear and leaves $\frac{1}{q} w^*$ invariant, 
we have $S_\tau\cN(\eta_0) = \frac{1}{q} w^* + S_\tau d$. 
By Lemma~\Ref{S:semig}, $\|S_\tau d\|_{2,\gamma-1} \le 
\e^{-(\gamma-3/2)\tau}\|d\|_{2,\gamma-1}$, so that $S_\tau d \to 0$ 
in $L^2_{\gamma-1}$ as $\tau \to \infty$. Thus, when $\tau$ is 
sufficiently large, we can apply Proposition~\ref{prop:C.NN} which 
gives 
\[
  \|\eta(\tau) - \eta^*_1\|_{2,\gamma-1} = \|\cN^{-1}({\ts\frac{1}{q}} 
  w^*{+}S_\tau d) - \eta^*_1\|_{2,\gamma-1} \le C_1\,\e^{-(\gamma-3/2)\tau},
\] 
for some $C_1 > 0$. This estimate holds in fact for all $\tau \ge 0$ 
with a possibly larger constant $C_1$, which proves \Ref{eq:prelim}.
Remark that, since $\eta(\tau)$ is nonnegative and $\|S_\tau d\|_1 
\to 0$, it follows from \Ref{eq:momid2} that $\| \eta(\tau)\|_{1,1} \to 
\frac{1}{\kappa} \,\e^{\euler}$ as $\tau \to \infty$. 

We next show that $\|\eta(\tau)\|_{2,\gamma}$ is uniformly bounded 
for all $\tau \ge 0$. We already know that $\|\eta(\tau)\|_{2,\gamma-1}$ 
and $\|\eta(\tau)\|_{1,1}$ remain bounded. Thus, using the integral 
equation \Ref{eq:Sinteq} together with the bounds \Ref{eq:semig} and 
\Ref{eq:bdd1}, we obtain
\[
  \|\eta(\tau)\|_{2,\gamma} \le \e^{-(\gamma-1/2)\tau} 
  \|\eta_0\|_{2,\gamma} + \int_0^\tau \beta(s) 
  \,\e^{-(\gamma-1/2)(\tau-s)} (q\|\eta(s)\|_{2,\gamma}+C)\dd s, 
\]
for some $C > 0$. Remark that $\beta(\tau) = \eta(\tau,1) = 
\cN(\eta(\tau))|_{y=1} = \frac{1}{q} + \e^\tau d(\e^\tau)$. Since 
$d \in L^2_{\gamma-1}$, it follows that $\beta(\tau) = \frac{1}{q} + 
\epsilon(\tau)$ with $\epsilon \in L^1(\R_+)$. Setting $H(\tau) = 
\e^{(\gamma-1/2)\tau} \|\eta(\tau)\|_{2,\gamma}$, we thus find
\[
  H(\tau) \le H(0) + \int_0^\tau (1{+}q|\epsilon(s)|)H(s)\dd s
  + C\,\e^{(\gamma-1/2)\tau} \FOR  \tau \ge 0, 
\]
for some $C > 0$. Since $\gamma-1/2 > 1$ and $\epsilon \in L^1(\R_+)$, 
it follows from Gronwall's lemma that $H(\tau) \le C_2 
\,\e^{(\gamma-1/2)\tau}$ for some $C_2 > 0$, hence 
$\|\eta(\tau)\|_{2,\gamma} \le C_2$ for all $\tau \ge 0$. 

Finally, if $\gamma \ge 2$, we show that $\eta(\tau)$ converges 
to $\eta^*_1$ in $L^2_\gamma$. To do this, we consider the integral
equation satisfied by $r(\tau) = \eta(\tau) - \eta^*_1$, namely
\[
  r(\tau) = S_\tau r(0) + \int_0^\tau S_{\tau-s} \Big\{\eps(s)
  T_1 \bbQ[\eta(s)] +\frac{1}{q} \big(T_1\bbQ[\eta^*_1{+}r(s)] - 
  T_1\bbQ[\eta^*_1]\big)\Big\}\dd s.
\]
In view of \Ref{eq:prelim} and Lemma~\ref{lem:bdd}, 
there exists $C_3 > 0$ such that
\[
  \|T_1 \bbQ[\eta(s)]\|_{2,\gamma} \le C_3, \quad
  \|T_1\bbQ[\eta^*_1{+}r(s)] - T_1\bbQ[\eta^*_1]\|_{2,\gamma} \le
   q \|r(s)\|_{2,\gamma} + C_3 \,\e^{-(\gamma-3/2)s}.
\]
Using Lemma~\ref{S:semig} again, we find that $R(\tau) = 
\|r(\tau)\|_{2,\gamma}$ satisfies the integral inequality
\[
  R(\tau) \le \e^{-(\gamma-1/2)\tau}R(0) + \int_0^\tau
  \e^{-(\gamma-1/2)(\tau-s)}\Big\{C_3|\eps(s)| + R(s) 
  + C_3 \e^{-(\gamma-3/2)s}\Big\}\dd s.
\]
Since $\eps(\tau)=\e^\tau d(\e^\tau)$ with $d\in L^2_{\gamma-1}$, 
we have
\[
  \int_0^\tau \e^{(\gamma-1/2)s}|\eps(s)|\dd s= \int_1^{\e^\tau}
  y^{\gamma-1/2}|d(y)| \dd y\leq \Big(\int_0^{\e^\tau} y\dd y
  \Big)^{1/2} \|d\|_{2,\gamma-1} \le \e^\tau \|d\|_{2,\gamma-1}, 
\]
hence there exists $C_4 > 0$ such that
\[
  R(\tau) \le C_4 \,\e^{-(\gamma-3/2)\tau} + \int_0^\tau
  \e^{-(\gamma-1/2)(\tau-s)} R(s) \dd s.
\]
Using Gronwall's lemma, we conclude that $R(\tau) \le C_5(1{+}\tau) 
\,\e^{-(\gamma-3/2)\tau}$ for some $C_5 > 0$, which is the desired 
result. 
\end{proof}

We now argue that the convergence towards the steady state $\eta^*_1$ 
cannot be faster than $\e^{-(\gamma-3/2)\tau}$ in the norm of $L^2_\gamma$, 
so that the result of Theorem~\ref{thm:conv} is optimal. To see 
this, we study the linearization of \Ref{eq:S.1} around $\eta^*_1$.
Setting $\eta(\tau) = \eta^*_1 + b(\tau)$, we obtain the linearized 
equation $\partial_\tau b = Ab$, where
\[
  (A b)(y)= (yb)'(y) + \Big(\frac{1}{q} T_1(\bbQ'[\eta^*_1]*b) 
  + b(1)T_1\bbQ[\eta^*_1]\Big)(y).
\]
Since we are interested in solutions $\eta(\tau) \in P_\gamma$, 
we study this operator in the space 
\[
  X_\gamma = \Bigset{b \in L^2_\gamma}{\int_1^\infty b(y)\dd y=0}.
\]

\begin{proposition}\label{prop:C.spec}
If $\gamma > 3/2$, the operator $A$ on $X_\gamma$ has
$\sigma = -(\gamma{-}3/2)$ in its spectrum.
\end{proposition}

\begin{proof}
For $\delta > \gamma{+}1/2$ we define a Lipschitz function 
$b_\delta \in X_\gamma$ by 
\[
  b_\delta(y)=\left\{\ba{cl} 
  -1 & \text{for } y \in [1,Y_\delta],\\
  -1 + (1{+}(Y_\delta{+}1)^{-\delta})(y-Y_\delta) & \text{for } 
    y \in (Y_\delta,Y_\delta+1),\\
   y^{-\delta}&\text{for } y \ge Y_\delta+1,\ea\right.
\]
where $Y_\delta$ is chosen such that $b_\delta$ has mean 0. Note that
$Y_\delta$ has a finite limit as $\delta\searrow \gamma{+}1/2$. 

Our aim is to show that  $A b_\delta + (\gamma{-}3/2)b_\delta$ stays
bounded in $X_\gamma$ as $\delta\searrow \gamma{+}1/2$, while 
$b_\delta$ is unbounded. For this purpose, we compute the asymptotic 
behavior of $A b_\delta$ as $y \to \infty$. Since $\eta^*_1$ decays 
faster than $\e^{-\lambda y}$ for some $\lambda > 0$ and since 
$\int_1^\infty \eta^*_1(y) \dd y = 1$, we obtain 
$(\bbQ'[\eta^*_1]*b_\delta)(y) = qb_\delta(y) + \cO(y^{-\delta-1})$ as 
$y\to \infty$, where $q=Q'(1)$. 
It follows that 
\[
\ba{rcl}
  (A b_\delta)(y) &=& (-\delta{+}1)y^{-\delta}+ (y{-}1)^{-\delta}
    + \cO(y^{-\delta-1}) + \cO(\e^{-\lambda y}) \\
  &=&(-\delta{+}2)y^{-\delta} + \cO(y^{-\delta-1}) \FOR  y \to \infty,
\ea
\]
where the remainder term is uniform in $\delta$ for $\delta \approx
\gamma{+}1/2$. This implies the estimate 
\[
  \|A b_\delta + (\gamma{-}3/2)b_\delta\|_{2,\gamma} \le
  (\delta{-}\gamma{-}1/2) \|b_\delta\|_{2,\gamma} + C \le 2C
\]
as $\delta\searrow \gamma{+}1/2$, since $\|b_\delta\|_{2,\gamma}
\approx 1/\sqrt{\delta{-}\gamma{-}1/2}$. This proves the claim. 
\end{proof}

To conclude this section, we also give a global stability result 
for the steady states $\eta^*_\theta$ with $0 < \theta < 1$. 

\begin{theorem}\label{thm:conv2}
Let $0 < \theta < 1$ and $\theta{+}1/2 < \gamma < \min\{3/2,2\theta{+}1/2\}$.
Assume that the initial value $\eta_0 \in \bbP$ satisfies 
$\eta_0 {-} \nu \eta^*_\theta \in L^2_{\gamma}$ for some $\nu > 0$, 
and let $\eta \in C^0([0,\infty),\bbP)$ be the solution of 
\Ref{eq:Sinteq} given by Theorem~\ref{S:Sexist}. Then there exists 
$C > 0$ such that 
\begin{equation}\label{eq:conv2}
  \|\eta(\tau)-\eta^*_\theta\|_{2,\gamma-\theta} \le C 
  \e^{-(\gamma -\theta -1/2)\tau} \FOR  \tau \ge 0. 
\end{equation}
\end{theorem}

\begin{remarks}\\
{\bf 1.} From \Ref{eq:asym1}, we know that $\eta^*_\theta(y) \sim 
y^{-1-\theta}$ as $y \to \infty$, so that $\eta^*_\theta \in 
L^2_{\gamma'}$ if and only if $\gamma' < \theta{+}1/2$. Thus, 
the assumption $\gamma > \theta{+}1/2$ guarantees that the 
difference $\eta_0 {-} \nu \eta^*_\theta$ decays faster than 
$\eta^*_\theta$ at infinity (otherwise, we could just choose $\eta_0 = 
\eta^*_{\theta'}$ for some $\theta' < \theta$, in which case 
$\eta(\tau) = \eta^*_{\theta'}$ for all $\tau \ge 0$ so that 
\Ref{eq:conv2} certainly fails.) For instance, if $\eta_0 \in 
\bbP \cap L^2$ is such that
\[
   \eta_0(y) = \frac{C}{y^{1+\theta}} + \cO\Bigl(\frac{1}
   {y^{1+\theta+\epsilon}}\Bigr) \FOR y \to \infty,
\]
where $C > 0$ and $\epsilon > 0$, then the assumptions of 
Theorem~\ref{thm:conv2} are satisfied for some $\nu$ and $\gamma$. 
On the other hand, the hypothesis $\gamma < 2\theta{+}1/2$ 
ensures that $\eta^*_\theta$ and hence $\eta_0$ lie in 
$L^2_{\gamma-\theta}$, so that $\eta(\tau) \in L^2_{\gamma-\theta}$
for all $\tau \ge 0$. 
\\[1mm]
{\bf 2.} Setting formally $\theta=1$ in \Ref{eq:conv2}, we recover 
\Ref{eq:prelim}. However, the main difference between the two results 
is the upper bound $\gamma < 3/2$ in Theorem~\ref{thm:conv2} 
which limits the decay rate in time of the perturbations. Even
for compactly supported perturbations, the convergence in 
\Ref{eq:conv2} is not faster than $\cO(\e^{-\delta\tau})$, 
where $\delta = \min\{\theta,1{-}\theta\}$. 
\end{remarks}

\begin{proof}
The proof is quite similar to that of \Ref{eq:prelim}, so we just 
indicate the main differences here. Proceeding as in the 
proof of Proposition~\ref{prop:C.N}, we first show that $\cN(\eta_0) = 
\frac{\theta}{q} w^*{+} d$ for some $d \in L^2_{\gamma-\theta}$. 
In analogy with \Ref{eq:repr}, we find
\[
  \wh d(\xi) = -\frac{1}{q} \log\Bigl(\frac{1-\Phi(\wh \eta_0(\xi))}
  {\e^{-\theta\wh w^*(\xi)}}\Bigr). 
\]
The crucial point is the behavior of $\wh d(\xi)$ as 
$\xi \to 0$, which we now analyze. By assumption, $\eta_0 =
\nu \eta^*_\theta + \zeta$ for some $\zeta \in L^2_{\gamma}$ with
$\int_1^\infty \zeta(y)\dd y = 1{-}\nu$. From \Ref{eq:etathdef}, 
we have
\[
  \wh \eta^*_\theta(\xi) = \Psi\bigl(1{-} \e^{-\theta \wh w^*(\xi)} 
  \bigr) = 1 - \e^{-\theta \wh w^*(\xi)} H\bigl(\e^{-\theta \wh w^*(\xi)}
  \bigr), 
\]
where $H : z \mapsto (1{-}\Psi(1{-}z))/z$ is analytic in a neighborhood
of zero, with $H(0) = \Psi'(1) = 1/\kappa$. We recall that 
$\e^{-\theta\wh w^*(\xi)} = (\i\xi)^\theta \e^{\theta\euler}
\e^{-\theta\chi(\i\xi)}$ where $\chi$ is entire, see \Ref{eq:E1exp}. 
Since $\gamma < 2\theta+1/2$, we deduce that $r_1 : \xi \mapsto 
H(\e^{-\theta \wh w^*(\xi)})$ belongs to $H^{\gamma-\theta}((-c,c))$
for some $c > 0$, and that $r_1(0) = 1/\kappa$. Next, we observe
that $\wh\zeta(\xi) = 1{-}\nu - r_2(\xi)$, where $r_2 \in 
H^\gamma((-c,c))$ and $r_2(0) = 0$. Using the fractional analog of 
Lemma~\ref{lem:Hardy}, we conclude that $r_2(\xi) = 
\e^{-\theta \wh w^*(\xi)} r_3(\xi)$, where $r_3 \in 
H^{\gamma-\theta}((-c,c))$ and $r_3(0) = 0$. (Here, we need 
$\gamma < 3/2$: if $\gamma > 3/2$, the claim would be false unless 
$r_2'(0) = 0$.) Summarizing, we have shown
\[
  \wh\eta_0(\xi) = \nu \wh\eta^*_\theta(\xi) + \wh\zeta(\xi)
  = 1 - \e^{-\theta \wh w^*(\xi)} \Bigl(\frac{\nu}{\kappa} + 
  r_4(\xi) \Bigr),
\]
where $r_4 \in H^{\gamma-\theta}((-c,c))$ and $r_4(0) = 0$.
We now apply the inverse map $\Phi = \Psi^{-1}$ which is analytic 
in a neighborhood of $1$ with $\Phi'(1) = \kappa$. Using the 
fact that $H^{\gamma-\theta}$ is an algebra, we obtain
\[
  \Phi(\wh\eta_0(\xi)) = 1 - \e^{-\theta \wh w^*(\xi)} 
  (\nu + r_5(\xi)),
\]
where $r_5$ has the same properties as $r_4$. Since $\wh d(\xi) 
= -\frac1q \log(\nu{+}r_5(\xi))$, we conclude that $\wh d \in 
H^{\gamma-\theta}((-c,c))$ with $\wh d(0) = -\frac1q\log\nu$. 
 
Now, from \Ref{eq:explicit} we have $\cN(\eta(\tau)) = 
S_\tau(\frac{\theta}{q} w^*{+}d) = \frac{\theta}{q} w^* + S_\tau d$,
and $\|S_\tau d\|_{2,\gamma-\theta} \le \e^{-(\gamma-\theta-1/2)\tau}
\|d\|_{2,\gamma-\theta}$ for all $\tau \ge 0$. Moreover, it is easy to 
check that Proposition~\ref{prop:C.NN} and its proof remain valid if 
we replace everywhere $w^*$ with $\theta w^*$, $\eta^*_1$ with 
$\eta^*_\theta$, and $\gamma$ with $\gamma' = \gamma-\theta < 
\theta{+}1/2$ (as is explained above, this inequality ensures that 
$\eta^*_\theta \in L^2_{\gamma'}$.) Thus, we conclude that 
$$
  \|\eta(\tau){-}\eta^*_\theta\|_{2,\gamma-\theta} = 
  \|\cN^{-1}(\ts{\frac{\theta}{q}}w^*{+}S_\tau d) 
  -\eta^*_\theta\|_{2,\gamma-\theta} \le C \|S_\tau 
  d\|_{2,\gamma-\theta} = \cO(\e^{-(\gamma-\theta-1/2)\tau}),
$$
as $\tau \to \infty$, which is the desired result.
\end{proof}


\appendix

\section{Bounds on the exponential integral}
\label{app:A}

Let $\wh w^*(\xi) = \Eone(\i\xi)$, where $\Eone(z) = \int_1^\infty 
y^{-1}\e^{-zy}\dd y$ is the exponential integral. The goal of this
section is to prove that
\begin{equation}\label{eq:A1}
  |1 - \e^{-\theta \wh w^*(\xi)}| < 1 \FOR \theta \in (0,1] \AND 
  \xi \in \bbL^- \setminus \{0\}.
\end{equation}
For $\theta = 1$ and $\xi \in \R$, this property is illustrated in 
Figure~\ref{figure2}.

\figurewithtex{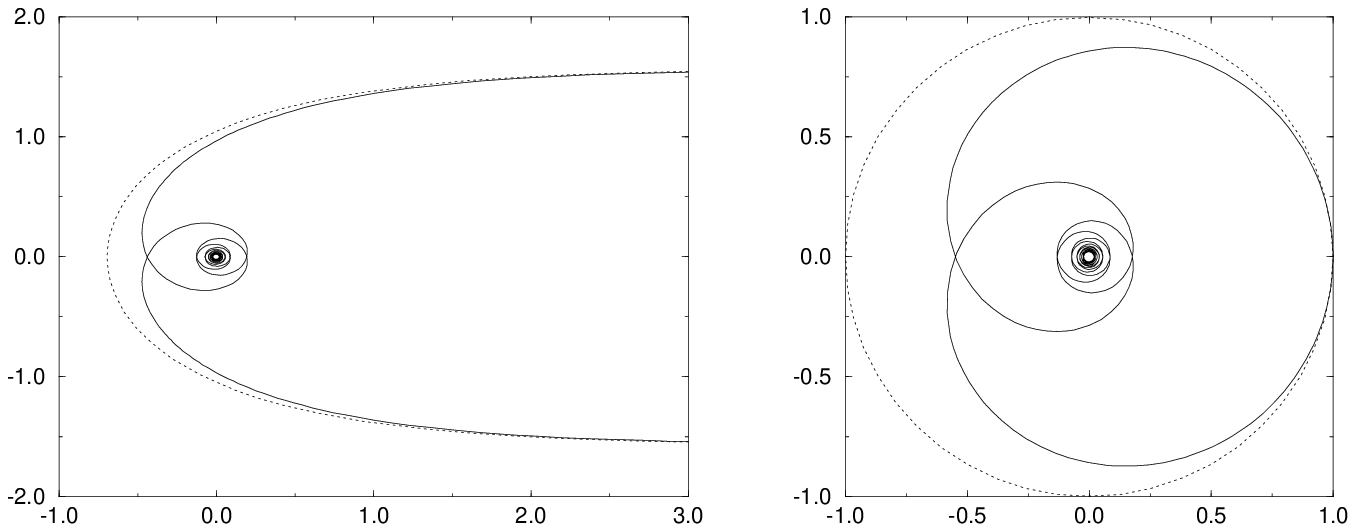}{fig2.tex}{6}{14}{figure2}{%
(a) The region  $D \subset \bbC$ delimited by the dotted line contains the
curve  $\set{\wh w^*(\xi) }{ \xi \in \bbR }$ (solid line). (b) The curve 
$\set{1{-}\e^{-\wh w^*(\xi)} }{ \xi \in \bbR }$ (solid line) is contained 
in the unit disk of $\bbC$.}

Fix $0 < \theta \le 1$, and define $F_\theta : \bbL^- \to \C$ by
$F_\theta(0) = 1$ and $F_\theta(\xi) = 1 - \e^{-\theta \wh w^*(\xi)}$
for $\xi \in \bbL^- \setminus \{0\}$. Then $F_\theta$ is continuous
on $\bbL^-$, and analytic in the interior of $\bbL^-$. Moreover, 
$F_\theta$ is uniformly bounded, because $|F_\theta(\xi)| \le 
1 + \exp(-\theta \Re(\wh w^*(\xi)))$ and 
\[
  \Re(\wh w^*(\xi)) = \int_1^\infty \frac{1}{y} \,\e^{\xi_2 y}
  \cos(\xi_1 y) \dd y \ge \int_{\pi/2}^\infty \frac{\cos(t)}{t}\dd t
  \approx - 0.472, 
\]
for all $\xi = \xi_1 {+} \i\xi_2$ with $\xi_1 \in \R$ and $\xi_2 < 0$. 
Finally, since $|\wh w^*(\xi)| \le \Eone(-\xi_2) \to 0$ as $\xi_2 \to 
-\infty$, it is clear that $|F_\theta(\xi)| \to 0$ as $\xi_2 \to
-\infty$, uniformly in $\xi_1 \in \R$. Thus, by the maximum modulus 
principle and the Phragmen-Lindel\"of theory (see e.g.\ \cite{Ru87}, 
Thm.~12.9), it is sufficient to show that \Ref{eq:A1} holds for 
all $\xi \in \R \setminus \{0\}$. 

Let $D \subset \C$ be the open region defined by
\[
  D = \set{x {+} \i y \in \C }{ |y| < \pi/2\,,~ x+\log(2\cos(y))  > 0 },
\]
see Figure~\ref{figure2}. 
As is easily verified, $w \in D$ implies $|1{-}\e^{-w}|
< 1$. Thus, all we need to show is that $\theta \wh w^*(\xi) \in D$ 
for all $\xi \in \R \setminus \{0\}$. Since $0 \in D$ and $D$ is 
convex, it is sufficient to prove this property for $\theta = 1$.

For $\xi > 0$ we define
\begin{equation}\label{eq:A2} 
  x(\xi) = \int_\xi^\infty \frac{\cos(t)}{t}\dd t \AND 
  y(\xi) = \int_\xi^\infty \frac{\sin(t)}{t}\dd t. 
\end{equation}
Then $\wh w^*(\xi) = \Eone(\i\xi) = x(\xi)-\i y(\xi)$ for $\xi > 0$ and
$\wh w^*(\xi) = x(|\xi|)+\i y(|\xi|)$ for $\xi < 0$. Moreover, 
$|y(\xi)| < \pi/2$ for all $\xi > 0$. Thus, it is enough to verify
that $K(\xi) > 0$ for all $\xi > 0$, where
\[
  K(\xi)  = x(\xi) + \log(2\cos(y(\xi))) \FOR  \xi > 0.
\]

We first observe that $K(\xi) > 0$ if $\xi > 0$ is sufficiently 
small. Indeed, in view of \Ref{eq:E1exp}, we have the expansions
\[
  x(\xi) = -\log\xi -\euler + \cO(\xi^2) \AND  
  y(\xi) = \frac{\pi}{2} -\xi + \cO(\xi^3) \FOR  \xi \searrow 0,
\]
hence $K(\xi) \to \log2 - \euler > 0$ as $\xi \searrow 0$. 

We next show that $K(\xi) > 0$ for $0 < \xi \le \pi/2$. If not, 
there would exist $\xi \in (0,\pi/2]$ such that $K(\xi) = 0$ and
$K'(\xi) \le 0$. In view of \Ref{eq:A2}, $K'(\xi) \le 0$
if and only if $\sin(\xi) \sin(y(\xi)) \le \cos(\xi) \cos(y(\xi))$. 
Since $0 < y(\xi) <\pi/2$, this is equivalent to $\xi + y(\xi) \le
\pi/2$, or $\cos(y(\xi)) \ge \sin(\xi)$. Therefore, $\xi \in
(0,\pi/2]$ should satisfy $x(\xi) + \log(2\sin(\xi)) \le K(\xi) = 0$. 
But this is impossible, because $x(\xi) + \log(2\sin(\xi)) \to 
\log 2-\euler > 0$ as $\xi \searrow 0$, and
\[
  \frac{\d}{\d \xi} \Bigl(x(\xi) + \log(2\sin(\xi))\Bigr) = 
  \frac{\xi - \sin(\xi)}{\xi \tan(\xi)} > 0 \quad \text{for } 
  0 < \xi < \pi/2. 
\]

It remains to show that $K(\xi) > 0$ for $\xi > \pi/2$. Let 
\[
  \bar x = -x(\pi/2) \approx 0.472 \AND  \bar y = \max\{y(\pi/2),
  -y(\pi)\} = -y(\pi) \approx 0.281. 
\]
(See \cite{AS72} for rigorous bounds on $x(\xi), y(\xi)$.) Using
the definitions \Ref{eq:A2}, it is easy to show that $|x(\xi)| \le \bar x$ 
and $|y(\xi)| \le \bar y$ for all $\xi \ge \pi/2$. Thus $x(\xi) + 
\log(2\cos(y(\xi))) \ge -\bar x + \log(2\cos(\bar y)) > 0$ for 
$\xi \ge \pi/2$. This concludes the proof. 


\section{Asymptotic behavior of the steady states}
\label{app:B}

Fix $\theta \in (0,1)$, and let $\eta = \eta^*_\theta : [1,\infty) \to
\R$ be the solution of \Ref{eq:time-indep} with $\beta = \theta/q$.
By Theorem~\ref{thm:equi}, $\eta$ is positive, strictly decreasing, and 
$\int_1^\infty \eta(y)\dd y = 1$. The aim of this section is to prove 
that the limit
\begin{equation}\label{eq:B1}
  L(\theta) = \lim_{y \to\infty} y^{1+\theta} \eta(y) 
\end{equation}
exists (and is finite). This is especially easy in the particular case
where $Q(z) = z$. Indeed, since $q=1$ and $\bbQ[\eta] = \eta$ in
this case, it follows from \Ref{eq:time-indep} that
\[
   0 = y\eta'(y) + \eta(y) + \theta \eta(y{-}1) \ge y\eta'(y) + 
   (1{+}\theta)\eta(y) \FOR  y \ge 2,
\]
hence $y \mapsto y^{1+\theta}\eta(y)$ is decreasing (and positive) 
for $y \ge 2$. In the general situation where $N = \deg Q > 1$, we
need the following estimate:

\begin{lemma}\label{lemB} For all $y \ge N$, 
\begin{equation}\label{eq:B2}
  \bbQ[\eta](y) \ge \eta(y) Q'\Bigl(\int_1^{y/N} \eta(x)\dd x\Bigr).
\end{equation}
\end{lemma}

\begin{proof}
The only property of $\eta$ that will be used in this proof is
that $\eta$ is nonnegative and non-increasing. Thus, by linearity
and monotonicity, it is sufficient to prove \Ref{eq:B2} in the case 
where $Q(z) = z^j$ for some $j \in \N$, $j \ge 2$. For $a \ge 1$, we 
denote
\begin{eqnarray*}
  D_j(a) &=& \set{(x_1,\dots,x_j) }{ 1 \le x_1, \dots, x_j \le a }
   = [1,a]^j, \\
  S_j(a) &=& \set{(x_1,\dots,x_j) }{ 1 \le x_1 \le  \dots \le x_j 
   \le a }. 
\end{eqnarray*}
Then, for $y \ge j$, we have
\begin{eqnarray*}
  \bbQ[\eta](y) &=& \int_{D_j(y+1-j)} \eta(x_1)\cdot{\dots}\cdot
    \eta(x_j) \delta(x_1+\dots+x_j-y) \dd^j x \\
  &=& j! \int_{S_j(y+1-j)} \eta(x_1)\cdot{\dots}\cdot
    \eta(x_j) \delta(x_1+\dots+x_j-y) \dd^j x,
\end{eqnarray*}
where $\delta$ denotes the Dirac measure. To obtain a lower bound, we 
replace $\eta(x_j)$ with $\eta(y)$ in the last integral, and we perform 
the (trivial) integration over $x_j$. We obtain
\[
  \bbQ[\eta](y) \ge (j!) \eta(y) \int_{R_{j-1}(y)}
  \eta(x_1)\cdot{\dots}\cdot \eta(x_{j-1}) \dd^{j-1} x, 
\]
where $R_{j-1}(y) = \set{(x_1,\dots,x_{j-1}) }{ (x_1,\dots,x_{j-1},
y{-}x_1{-}\dots{-}x_{j-1}) \in S_j(y{+}1{-}j) }$. Now, it is straightforward 
to verify that $R_{j-1}(y) \supset S_{j-1}(y/j)$. Thus
\begin{eqnarray*}
  \bbQ[\eta](y) &\ge& (j!) \eta(y) \int_{S_{j-1}(y/j)} 
   \eta(x_1)\cdot{\dots}\cdot\eta(x_{j-1}) \dd^{j-1} x \\
  &=& j \eta(y) \int_{D_{j-1}(y/j)} 
   \eta(x_1)\cdot{\dots}\cdot\eta(x_{j-1}) \dd^{j-1} x \\
  &=& j \eta(y) \Bigl(\int_1^{y/j} \eta(x) \dd x \Bigr)^{j-1}
  \,=\, \eta(y) Q'\Bigl(\int_1^{y/j} \eta(x) \dd x \Bigr).
\end{eqnarray*}
This concludes the proof.
\end{proof}

Combining \Ref{eq:time-indep} and Lemma~\ref{lemB}, we obtain the
inequality
\[
   y\eta'(y) + \eta(y) + \frac{\theta}{q} \eta(y{-}1)\,
   Q'\Bigl(\int_1^\frac{y{-}1}{N} \eta(x)\dd x\Bigr) \le 0 
  \FOR y \ge N{+}1, 
\]
where $\eta(y{-}1)$ may also be replaced by $\eta(y)$. It follows that
\[
  \frac{\d}{\d y}(y^{1+\theta/2}\eta(y)) \le \theta y^{\theta/2} 
  \eta(y) \Bigl(\frac{1}{2} - \frac{1}{q} Q'\Bigl(\int_1^{\frac{y-1}{N}} 
  \eta(x)\dd x\Bigr)\Bigr) \FOR  y \ge N{+}1. 
\]
Since $\int_1^\infty \eta(y)\dd y = 1$ and $q = Q'(1)$, the right-hand 
side becomes negative for $y$ sufficiently large. This shows that
\begin{equation}\label{eq:previous}
  \sup_{y \ge 1} y^{1+\theta/2}\eta(y) < \infty. 
\end{equation}
Similarly, for $y \ge N{+}1$,
\begin{equation}\label{eq:B3}
  \frac{\d}{\d y}(y^{1+\theta}\eta(y)) \le  y^{1+\theta} 
  \eta(y) \Theta(y)\quad \text{with }
\Theta(y)=\frac{\theta}{y}\Big( 1 - \frac{1}{q} Q'\Bigl(\int_1^{\frac{y-1}{N}} 
  \eta(x)\dd x\Bigr)\Big)\geq 0.
\end{equation}
It follows from \Ref{eq:previous} that $1-\int_1^{(y-1)/N}\eta(x)\dd
 x=\cO(y^{-\theta/2}) $ for $y \to \infty$, which yields
 $\Theta(y) =
\cO(y^{-1-\theta/2})$  and hence $\Theta\in L^1((N{+}1,\infty))$.  
Thus, the differential inequality 
\Ref{eq:B3} implies that the limit \Ref{eq:B1} exists. 


\section{Bounds on the nonlinearity}
\label{app:C}

In this section, we sketch the proofs of Lemmas~\ref{S:nonlin} 
and \ref{lem:bdd}. Without loss of generality, we assume here
that $Q(z) = z^m$ for some $m \in \N_*$ (the general case follows
by linearity). To bound the convolution products, we repeatedly 
use Young's inequality $\|f*g\|_p \le \|f\|_p \|g\|_1$ where 
$f \in L^p$ and $g \in L^1$. 

\smallskip
\noindent{\bf Proof of Lemma~\ref{S:nonlin}.}\\
{\bf a)} If $\eta \in L^1$, then $\bbQ[\eta] = 
\eta^{*m} \in L^1$ and $\|\bbQ[\eta]\|_1 \le \|\eta\|_1^m = 
Q(\|\eta\|_1)$. If $\eta,\tilde\eta \in L^1$, then
\begin{equation}\label{eq:diffeta}
  \bbQ[\eta]-\bbQ[\tilde\eta] \,=\, (\eta{-}\tilde\eta) * 
  \eta * \dots *\eta \,+\, \dots \,+\, \tilde \eta * \tilde 
  \eta * \dots * (\eta{-}\tilde\eta) 
\end{equation}
($m$ terms of $m$ factors), hence $\|\bbQ[\eta]-\bbQ[\tilde\eta]\|_1 
\le m r^{m-1} \|\eta-\tilde\eta\|_1 = Q'(r)\|\eta-\tilde\eta\|_1$.

\noindent{\bf b)} Assume now that $\eta \in L^p_\gamma \hookrightarrow
L^1$. For all $y \ge 1$, 
\[
  y^\gamma (T_1\bbQ[\eta])(y) = \int_{\R^m} 
  \eta(x_1)\dots\eta(x_m) (1{+}x_1{+}\dots{+}x_m)^\gamma 
  \,\delta(1{+}x_1{+}\dots{+}x_m{-}y) \dd^m x,
\]
where $\delta$ denotes the Dirac measure. Due to the support property
of $\eta$, only the values $x_1, \dots, x_m \ge 1$ contribute to the 
integral. For such values, we have the estimate
\begin{equation}\label{eq:crude}
  (1{+}x_1{+}\dots{+}x_m)^\gamma \le C(x_1^\gamma {+} \dots {+}
  x_m^\gamma), 
\end{equation}
where $C > 0$ depends on $m,\gamma$. Thus, $|y^\gamma(T_1\bbQ[\eta])|$
is bounded by a sum of $m$ convolution products of the form 
$|y^\gamma \eta| * |\eta|^{*(m-1)}$. Taking the $L^p$ norm and
using Young's inequality, we obtain
\[
  \|T_1\bbQ[\eta]\|_{p,\gamma} \le C m \|\eta\|_1^{m-1} 
  \|\eta\|_{p,\gamma} = C Q'(\|\eta\|_1) \|\eta\|_{p,\gamma}.
\]
Finally, using the decomposition \Ref{eq:diffeta} and proceeding as
above, we find
\begin{eqnarray*}
 \|T_1\bbQ[\eta] - T_1\bbQ[\tilde\eta]\|_{p,\gamma} 
 & \le & C(m r^{m-1}\|\eta-\tilde\eta\|_{p,\gamma}
 + m(m{-}1)r^{m-2} R \|\eta-\tilde\eta\|_1 ) \\
 & \le & C(Q'(r) + R Q''(r)) \|\eta-\tilde\eta\|_{p,\gamma},
\end{eqnarray*}
where $R = \max\{\|\eta\|_{p,\gamma},\|\tilde\eta\|_{p,\gamma} \}$
and $r = \max\{ \|\eta\|_1,\|\tilde\eta\|_1\}$. Since $r \leq R$ and
$RQ''(R)\leq CQ'(R)$, this is the desired result.
\QED

\smallskip
\noindent{\bf Proof of Lemma~\ref{lem:bdd}.} 
The proof follows the same lines, except that \Ref{eq:crude} is
replaced with a different estimate, which can be established 
by induction over $m$. If $\gamma \ge 1$ and $m \in \N_*$, there 
exists $C > 0$ such that, for all $x_1, \dots,x_m \ge 1$,
\[
   (1{+}x_1{+}\dots{+}x_m)^\gamma \,\le\, \sum_{i=1}^m \Bigl( 
   x_i^\gamma + C x_i^{\gamma-1}\prod_{j\neq i} x_j\Bigr). 
\]
{\bf a)} If $\gamma > 3/2$ and $\eta \in P_\gamma$, then
\begin{eqnarray*}
 \|T_1\bbQ[\eta]\|_{2,\gamma} 
  & \le & m \|\eta\|_1^{m-1} \|\eta\|_{2,\gamma} +
   C m \|\eta\|_{1,1}^{m-1} \|\eta\|_{2,\gamma-1}\\  
  & = & Q'(1) \|\eta\|_{2,\gamma} +
   C Q'(\|\eta\|_{1,1}^{m-1}) \|\eta\|_{2,\gamma-1}.\\  
\end{eqnarray*}
{\bf b)} If $\eta,\tilde\eta \in P_\gamma$, let 
$M = \max\{ \|\eta\|_{2,\gamma},\|\tilde\eta\|_{2,\gamma}\}$, 
$r_1 = \max\{\|\eta\|_{1,1},\|\tilde\eta\|_{1,1}\} \le M$, and
$r = \max\{\|\eta\|_1,\|\tilde\eta\|_1\} = 1$. Then 
\begin{eqnarray*}
 \|T_1\bbQ[\eta] - T_1\bbQ[\tilde\eta]\|_{2,\gamma} 
 & \le & m r^{m-1} \|\eta-\tilde\eta\|_{2,\gamma}
   + m(m{-}1)r^{m-2} M \|\eta-\tilde\eta\|_1  \\
 && +\, Cm r_1^{m-1} \|\eta-\tilde\eta\|_{2,\gamma-1}
   + Cm(m{-}1)r_1^{m-2} M \|\eta-\tilde\eta\|_{2,1}  \\
 & = & Q'(1)\|\eta-\tilde\eta\|_{2,\gamma}
   + Q''(1) M \|\eta-\tilde\eta\|_1  \\
 && +\, C Q'(r_1) \|\eta-\tilde\eta\|_{2,\gamma-1}
   + C Q''(r_1) M \|\eta-\tilde\eta\|_{2,1} .
\end{eqnarray*}
If $\gamma \ge 2$, the last three terms in the right-hand side
can be bounded by $CQ'(M)\|\eta-\tilde\eta\|_{2,\gamma-1}$. \QED




\parskip 0pt
\parindent 30pt

\begin{thebibliography}{NK99}

\bibitem[AS72]{AS72}
\textsc{M. Abramowitz, I. Stegun:} \textit{Handbook of Mathematical
Functions.} Dover, 1972. 

\bibitem[BDG94]{BDG94} 
\textsc{A. Bray, B. Derrida, C. Godr\`eche:} Non-trivial algebraic decay
in a soluble model of coarsening. \textit{Europhys. Letters}
\textbf{27} (1994) 175--180.

\bibitem[BrD95]{BrD95}
\textsc{A. Bray, B. Derrida:} Exact exponent $\lambda$ of the 
autocorrelation function for a soluble model of coarsening. 
\textit{Phys. Rev.} \textbf{E 51} (1995) 1633--1636.  

\bibitem[CaP89]{CaP89}
\textsc{J. Carr, R.L. Pego:} Metastable patterns in solutions 
of $u_t = \epsilon^2 u_{xx}-f(u)$. {\it Commun. Pure Appl. Math.} 
\textbf{42} (1989) 523--576. 

\bibitem[CaP92]{CaP92}
\textsc{J. Carr, R.L. Pego:} Self-similarity in a coarsening model in one
  dimension.  {\it Proc. Roy. Soc. Lond.} \textbf{A 436} (1992) 569--583. 

\bibitem[CaP00]{CaP00}
\textsc{J. Carr, R.L. Pego:} Self-similarity in a cut-and-paste model of
  coarsening, {\it Proc. Roy. Soc. Lond.} \textbf{A 456} (2000)
1281--1290.  

\bibitem[De97]{De97} 
\textsc{B. Derrida:} Non-trivial exponents in coarsening phenomena.
\textit{Physica D} \textbf{103} (1997) 466--477. 

\bibitem[DGY91]{DGY91}
\textsc{B. Derrida, C. Godr\`eche, I. Yekutieli:} Scale-invariant
regimes in one-dimensional models of growing and coalescing droplets.
\textit{Phys. Rev.} \textbf{A 44} (1991) 6241--6251. 

\bibitem[Do50]{Do50}
\textsc{G. Doetsch:} \textit{Handbuch der Laplace-Transformation, 
Band I.} Birkh\"auser, Basel, 1950.  

\bibitem[Fe71]{Fe71}
\textsc{W. Feller:} \textit{An Introduction to Probability Theory and
its Application, Vol. II}. Second edition, Wiley, New-York, 1971.    

\bibitem[HLP59]{HLP59}
\textsc{G. Hardy, J. Littlewood, G. P\'olya:} \textit{Inequalities}. 
Third edition, Cambridge University Press, 1959.    

\bibitem[KBN97]{KBN97}
\textsc{P.L. Krapivsky, E. Ben-Naim:} Domain statistics in coarsening
systems. \textit{Phys. Rev. E} \textbf{56} (1997) 3788--3798.  

\bibitem[LiS61]{LiS61}
\textsc{I.M. Lifshitz, V.V. Slyozov:} The kinetics of precipitation 
from supersaturated solid solutions. \textit{J. Phys. Chem. Solids} 
\textbf{19} (1961) 35--50.  

\bibitem[NaK86]{NaK86}
\textsc{T. Nagai, K. Kawasaki:}  Statistical
dynamics of interacting kinks. II \textit{Physica} \textbf{A 134}
(1986) 483--521.  

\bibitem[PeR92]{PeR92}
\textsc{K. Pesz, G.J. Rodgers:} Kinetics of growing and coalescing
droplets. \textit{J. Phys. A} \textbf{25} (1992) 705--713.  

\bibitem[Ru87]{Ru87} \textsc{W. Rudin:} \textit{Real and Complex
Analysis}. Third edition, McGraw-Hill, New York, 1987. 

\bibitem[RuB94]{RuB94}
\textsc{A. Rutenberg, A. Bray:} Phase-ordering kinetics of
one-dimensional nonconserved scalar systems. \textit{Phys. Rev.} 
\textbf{E 50} (1994) 1900--1911.   

\bibitem[Sch66]{Sch66} \textsc{L. Schwartz:} \textit{Th\'eorie des 
Distributions}. Second edition, Hermann, Paris, 1966. 

\bibitem[Vo85]{Vo85}
\textsc{P. Voorhees:} The theory of Ostwald ripening. \textit{J. Stat.
Phys.} \textbf{38} (1985) 231--252.  

\end{thebibliography}
\end{document}